\documentclass[11pt]{article}
\usepackage[T1]{fontenc}
\usepackage[utf8]{inputenc}
\usepackage[a4paper, total={5.5in, 9in}]{geometry}

\usepackage{siunitx}              
\usepackage{hyperref}

\usepackage{graphicx}%
\usepackage{amsmath,amssymb,amsfonts}%
\usepackage[capitalize,nameinlink]{cleveref}
\usepackage{amsthm}%
\usepackage{mathrsfs}%
\usepackage[title]{appendix}%
\usepackage{xcolor}%
\usepackage{textcomp}%
\usepackage{manyfoot}%
\usepackage{booktabs}%
\usepackage{algorithm}%
\usepackage{algorithmicx}%
\usepackage{algpseudocode}%
\usepackage{listings}%
\usepackage[normalem]{ulem}
\usepackage{subcaption} 
\usepackage{booktabs}
\usepackage{multirow}
\usepackage{mathtools}
\usepackage{microtype}
\usepackage{makecell}

\newtheorem{definition}{Definition}%

\usepackage[nocitation]{apacite}
\usepackage{natbib}
\usepackage{authblk} 

\providecommand{\keywords}[1]{\textbf{\textit{Keywords---}} #1}

\title{On leveraging constrained smooth additive regression models for global optimization}
\author[1]{Marina Cuesta}
\author[2]{Claudia D'Ambrosio}
\author[1]{María Durban}
\author[1]{Vanesa Guerrero\thanks{Corresponding author: \texttt{vanesa.guerrero@uc3m.es}}}
\author[2]{Renan Spencer Trindade}

\affil[1]{Department of Statistics, Universidad Carlos III de Madrid, Calle Madrid 126, 28901 Getafe, Spain}
\affil[2]{LIX CNRS, École Polytechnique, Institut Polytechnique de Paris, Route de Saclay, Palaiseau, 91128, France}

\date{}

\begin{document}

\maketitle

\begin{abstract}
 Many real-world decision-making processes rely on solving mixed-integer nonlinear programs (MINLPs). However,  finding high-quality solutions to MINLPs is often computationally demanding, motivating the development of specialized algorithms to improve their tractability. In this work, we propose Mixed-Integer Smoothing Surrogate Optimization with Constraints (MISSOC), a novel optimization algorithm that builds and solves  approximations of challenging MINLPs. MISSOC approximates complicating functions in an MINLP using smooth additive regression models with \unboldmath{$B-$}splines. Expert knowledge can be incorporated into the approximating functions through shape constraints related to bounds, monotonicity and curvature over the observed domain. A surrogate of the original problem is then obtained by replacing the original complicating functions with their approximations, making it more tractable in practice. MISSOC presents an innovative integration of statistical modeling into mathematical optimization and fills a gap in the literature by building surrogates that are both data-driven and knowledge-driven.
The proposed algorithm is illustrated on the real-world Water Distribution Network problem and evaluated through a set of experiments that include benchmark instances and the real-world Hydro Unit Commitment problem. Together, they demonstrate that MISSOC handles MINLPs with integer variables and complicating functions appearing in the objective or in the constraints. MISSOC is  evaluated with different state-of-the-art solvers and with the Sequential Convex MINLP (SC-MINLP) algorithm. The latter exploits the separable structure of the approximating functions, which are sums of piecewise univariate polynomials. The experiments show that MISSOC can obtain high-quality solutions for challenging MINLPs, particularly when used in combination with the SC-MINLP algorithm.

\vspace{0.4 cm}
\noindent \keywords{Mixed-Integer nonlinear program, Surrogate functions, Shape-constrained regression}

\end{abstract}

\section{Introduction}\label{sec:intro}

Mixed-Integer nonlinear programs (MINLPs) are a powerful and versatile mathematical optimization framework that enables the modeling and solving of challenging  real-world problems involving discrete and continuous decisions, as well as nonlinear relationships.  As such, it plays a key role in a wide range of applications, including engineering, finance, and many other fields \citep{borghetti2015optimal,bussieck2003minlplib,Taktak2017}. However, finding a global solution to a general MINLP, even in the case of finite bounds on the variables, is NP-hard \citep{burer2012}, making these problems computationally challenging. Global optimization solvers usually succeed in optimally solving small instances of nonconvex MINLPs, namely MINLPs having a nonconvex continuous relaxation, obtained by relaxing integrality requirements. However, they may struggle with medium- to large-size instances. This fact motivates the development of specialized methods and reformulations for subclasses of MINLPs that aim to improve tractability and extend the range of solvable problem instances. In particular, one of these approaches consists in obtaining an approximation of the original MINLP that is more tractable to solve. This surrogate  problem can be obtained in a data-driven (e.g., \cite{bertsimas2025global,bertsimas2023global}) or knowledge-driven (e.g., \cite{Codsi2025,duguet-ngueveu2022,warwicker2024efficient}) way.
In this paper, we develop a new framework for building and solving surrogate problems for MINLPs that are both data- and knowledge-driven.
In the rest of the manuscript the concept of {\it approximation} and {\it surrogate} will be used as synonyms.

Some approaches in the literature construct surrogate problems for MINLPs using data-driven techniques, typically based on machine learning models that can be straightforwardly embedded into optimization frameworks. The integration of the objective function or constraints derived from data-driven methods has been investigated in recent works such as \cite{bertsimas2023global} and \cite{bertsimas2025global}. In these studies, nonlinear constraints are approximated using sampling-based techniques to train a machine learning model. The resulting approximations are incorporated into the mathematical optimization model to yield a mixed-integer linear  program (MILP) surrogate of the original MINLP. In \cite{bertsimas2023global}, hyperplane-based decision-trees are used to derive the approximating functions, while in \cite{bertsimas2025global} the authors consider gradient boosted trees, multi layer perceptrons, and support vector machines. The resulting surrogate  problem is formulated as an MILP and is then solved. Then,  its optimal solution is used as a starting point for a local repair phase, where the potential infeasibility and suboptimality of the solution are either fixed or improved. However, these approaches have several limitations, including  that the authors focus their studies to MILP-representable approximations, which may be too simplistic to capture the structure of the original MINLP \citep{geissler2011using}. Moreover, they build the approximations using Interpretable AI\footnote{\url{https://www.interpretable.ai/}}, a commercial software. Thus, it might be difficult for a standard user to integrate expert knowledge in the machine learning training phase.

Several different knowledge-driven approaches for building approximations of MINLPs have been proposed in the literature. For instance, \cite{Rebennack2015} focus on nonlinear functions of two variables. They develop piecewise linear (PWL) approximations that ensure the resulting function either under- or over-estimates the original bivariate function. Additionally, they guarantee that the approximation error remains within a given tolerance. \cite{Codsi2025} focus on univariate nonlinear functions. They propose discontinuous PWL approximations that ensure a bounded, pointwise approximation error while minimizing the number of pieces of the PWL approximation. Their idea is that, in standard MILP formulations of PWL approximations, the number of binary variables is proportional to the number of pieces. Thus, by minimizing the number of pieces, they expect to have a more tractable MILP approximation. They validate their approach on a testbed of univariate nonlinear functions and on a class of MINLPs with linear constraints. \cite{duguet-ngueveu2022} extend the methodology by \cite{Codsi2025} (previously available as the technical report \cite{codsi2021lina}) to bivariate functions. Their approach guarantees that the domain of the nonlinear function is split into the minimum number of pieces required to ensure a bounded approximation error. The reader is referred to \cite{geissler2011using} for more details on approaches based on PWL approximation of MINLPs. The recent work by \cite{goss2026parabolic} proposes quadratic approximations for MINLPs, specifically parabolic approximations and relaxations. To derive these approximations, in the fitting phase, the authors solve MILPs. Their two-phase approach is evaluated on the MINLPlib benchmark \citep{bussieck2003minlplib}. The main limitations of the mentioned knowledge-driven approaches are as follows. First,  the aforementioned PWL approaches are tailored for specific classes of nonlinear functions, namely univariate or bivariate nonlinear functions. Second, the parabolic approximation method by \cite{goss2026parabolic} shows a time consuming phase to identify such an approximation. In particular, in their case it requires solving an MILP. Thus, in practice they precompute parabolic approximations of some recurring nonlinear functions and store them in lookup tables.

This paper aims to fill a gap in the literature by developing surrogate problems for MINLPs that are both data-driven \textit{and} knowledge-driven. In other words, given an MINLP, our goal is to solve it by approximating its complicating functions with simpler ones, and then obtaining a new MINLP, i.e. the \emph{surrogate MINLP}, which is  more tractable. Our aim is that this surrogate MINLP remains accurate while also allowing the incorporation of expert knowledge about the shape or other properties of the original functions being approximated. The tractability of the so-obtained surrogate MINLP is crucial to get high-quality solutions for large instances of MINLPs \citep{bhosekar2018advances}. In particular, we use constrained smooth additive regression models to gain leverage in tractability.

In this work, we develop a novel approach for incorporating shape constraints, representing expert knowledge, into the estimated functions that approximate the complicating functions in the original MINLP, while still accurately fitting the training data sampled from those functions. To do so, the works by \cite{Navarro2023AMC,Navarro2024multispline} are extended to the case of constrained smooth additive regression models. $B-$spline  functions \citep{de1978practical} are used to represent the univariate functions in the additive model and the smoothness of the curve is controlled through the number and location of the knots needed to define them. Functions estimated this way are sums of univariate piecewise polynomials that can be embedded into a mathematical programming model by means of a multiple choice formulation \citep{vielma2009mixed}. The number of knots translates into the number of binary variables needed to model them, while their additive and separable nature allows us to use customized solution approaches, such as the Sequential Convex-MINLP (SC-MINLP) algorithm \citep{d2012algorithmic,d2019strengthening, trindade2022comparing}. Building on the above, we propose the Mixed-Integer Smoothing Surrogate Optimization with Constraints (MISSOC) algorithm, which builds and solves surrogate MINLPs by approximating challenging functions using shape-constrained smooth additive regression models. MISSOC is more flexible than existing approaches in the literature as it preserves nonlinearities in the surrogate MINLP and combines both data-driven and knowledge-driven elements. This flexibility, however, comes at a cost: while the knowledge-driven methods discussed above build approximations of MINLPs with a prescribed maximum error bound over the full domain, our approach relies only on sampled data and the least squares method to estimate the approximating function. As a result, the approximation error is controlled only on the sampled points through the least squares criterion and we cannot guarantee a maximum error bound over the entire domain.

The remainder of this paper is organized as follows. Section~\ref{sec:preliminary_concepts} outlines the basic definitions of smooth regression models. We present our novel approach to estimate smooth additive regression models with constraints in Section~\ref{sec:shape_constrained_regression}. Then, the MISSOC approach to solve MINLPs is explained in detail in Section~\ref{sec:surrog}. Section~\ref{sec:demonstrative} presents a demonstrative example of MISSOC using the real-world Water Distribution Network problem \citep{bragalli2012optimal}. The proposed methodology is then tested in Section~\ref{sec:comput_results} on MINLPlib instances \citep{bussieck2003minlplib}, as well as in the real-world  Hydro Unit Commitment problem \citep{borghetti2015optimal,Taktak2017}. Finally, Section~\ref{sec:conclusions}  concludes the paper with final remarks and perspectives.

\section{Preliminary concepts}
\label{sec:preliminary_concepts}

This section introduces the basic concepts about smooth regression models that are needed to develop our new approach for smooth additive regression models with constraints in Section~\ref{sec:shape_constrained_regression}. Recall that this modeling approach is used to approximate complicating functions in MINLPs and to obtain data- and knowledge-driven  surrogate MINLPs. First, in Section~\ref{sec:preliminaries:univariate_regression} we review how to estimate smooth univariate regression models with $B-$splines. Then, Section~\ref{sec:preliminaries:univariate_shape_constraints} is devoted to review the approach by \cite{Navarro2023AMC} to include shape constraints, i.e. about the sign, monotonicity or curvature of the so-obtained curves. Finally, smooth additive regression models with $B-$splines are introduced in Section~\ref{sec:preliminaries:smooth_additive_regresion}.

\subsection{Smooth univariate regression models}
\label{sec:preliminaries:univariate_regression}

Let us consider a set of $n$ observations $\{(x_i, y_i) \in \mathbb{R}^2, \ i = 1, \dots, n \}$ drawn from a continuous covariate $X$  (the independent variable serving as input) and a continuous response variable $Y$, which represents the dependent output value we aim to predict, respectively. For simplicity, it can be assumed that $x_1 \leq x_2 \leq \dots \leq x_n$ without loss of generality. A univariate smooth regression model is given by

    \begin{equation}
        y_i = \alpha + f(x_i) + \epsilon_i, \quad i = 1, \dots, n,
       \label{eq:preliminaries:univariate_model}
    \end{equation}

\noindent where $\alpha \in \mathbb{R}$ is the intercept of the model,   $f: [x_{1}, x_{n}] \subset \mathbb{R} \rightarrow \mathbb{R}$ is a smooth function and $\epsilon_i \in \mathbb{R}$ are independent error terms with zero mean ($i =1,\dots,n$). 
 The error term $\epsilon_i$ in~(1) accounts for observational noise in the regression model. 

The aim is to accurately estimate $\alpha$ and $f$ in~\eqref{eq:preliminaries:univariate_model} based on the sampled values of $X$ and $Y$. This is a nontrivial task and many methods have been proposed in the literature \citep{schimek2013smoothing}. Among them, this work adopts a $B-$spline approach \citep{de1978practical},  since they provide a convenient representation of smooth functions through locally supported piecewise polynomials. Moreover, $B-$splines are low-rank smoothers, in the sense that they use a basis of significantly lower dimension than the number of data points, which provides regularization and computational efficiency. $B-$splines can be defined in multiple dimensions, but we focus on the univariate case. 

\begin{figure}[h!]
    \centering
    \begin{subfigure}{0.32\textwidth}
        \includegraphics[width=0.8\linewidth]{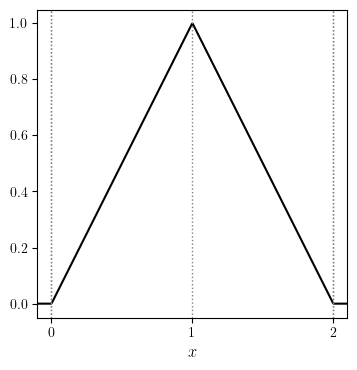}
        \caption{$d=1$. \label{subfig:d1}}
    \end{subfigure}
    \hfill
    \begin{subfigure}{0.32\textwidth}
        \includegraphics[width=0.8\linewidth]{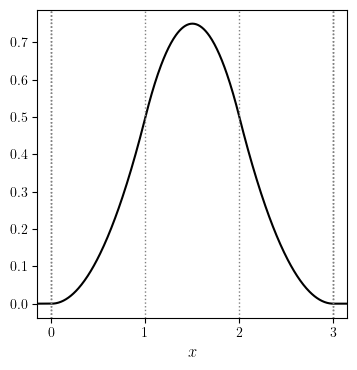}
        \caption{$d=2$. \label{subfig:d2}}
    \end{subfigure}
    \hfill
    \begin{subfigure}{0.32\textwidth}
        \includegraphics[width=0.8\linewidth]{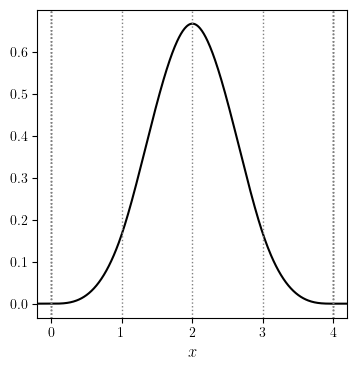}
        \caption{$d=3$. \label{subfig:d3}}
    \end{subfigure}
    \caption{$B-$spline functions of degrees 1, 2 and 3. The positions of the knots are indicated with vertical dotted lines.}
    \label{fig:preliminaries:Bsplines}
\end{figure}

 $B-$splines are defined with respect to a nondecreasing sequence of real values, called \emph{knots}, which partition the domain into intervals. A $B-$spline function of degree $d$ is a univariate piecewise polynomial that is nonzero only over the span of $d + 2$ consecutive knots where it consists of $d+1$ polynomial segments of degree $d$, joined at $d$ interior knots. At these points, the derivatives up to order $d - 1$ are continuous. Figure~\ref{fig:preliminaries:Bsplines} illustrates $B-$spline functions of degrees 1, 2 and 3.

 Given a polynomial degree $d$ and a nondecreasing knot sequence $t=\{t_1,\ldots,t_m\}$, the associated $B-$spline functions are uniquely determined via the Cox--de~Boor  \citep{de1978practical} recursive construction. These functions constitute the elementary building blocks used in the remainder of this section. In particular, for a degree $d$ and a knot sequence $t$, the $B-$spline basis consists of $m-d-1$ functions. 
For $l=1,\ldots,m-d-1$, the $l$-th $B-$spline function is defined recursively by

\begin{equation*}
B_{l,d,\mathbf{t}}(x) = \frac{x - t_{l}}{t_{l+d} - t_{l}} B_{l,d-1,\mathbf{t}}(x) + \frac{t_{l+d+1} - x}{t_{l+d+1} - t_{l+1}} B_{l+1,d-1,\mathbf{t}}(x),
\label{eq:preliminaries:B-spline_basis}
\end{equation*}

\noindent where
\begin{equation*}
B_{l,0,\mathbf{t}}(x) =
\begin{cases}
1 & \text{if } x \in [t_{l}, t_{l+1}), \\
0 & \text{otherwise}.
\end{cases}
\end{equation*}

\begin{figure}[h!]
    \centering
    \begin{subfigure}[t]{0.48\textwidth}
        \includegraphics[width=\linewidth]{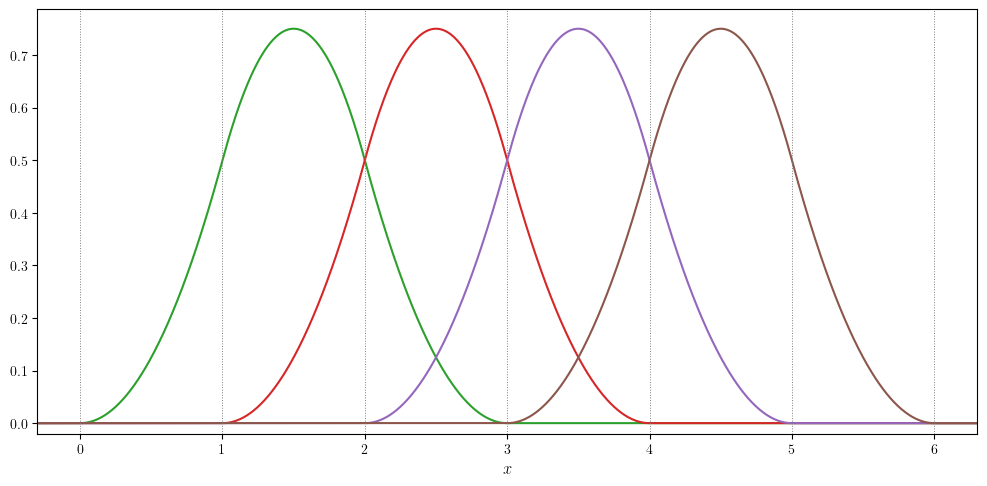}
        \caption{$B-$spline basis of degree  $d=2$. The positions of the knots are indicated with vertical dotted lines.}
        \label{fig:preliminaries:Bspline_basis}
    \end{subfigure}    
    \hfill
    \begin{subfigure}[t]{0.48\textwidth}
        \includegraphics[width=\linewidth]{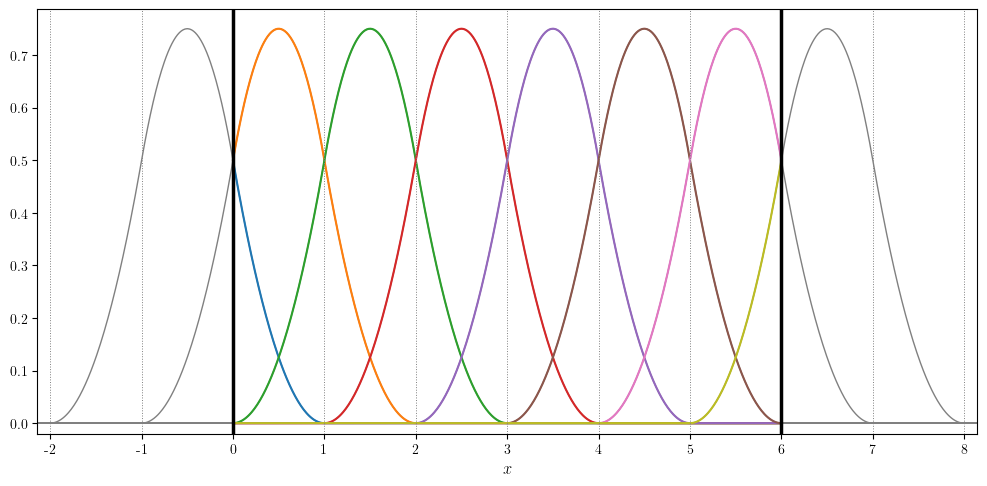}
        \caption{$B-$spline basis of degree $d=~2$ with extended knot sequence. Internal knot boundaries are highlighted with thick black vertical lines. The segments of the $B-$spline functions lying in the external intervals are colored in grey.}
        \label{fig:preliminaries:Bspline_basis_extended}
    \end{subfigure}    
    \caption{Example of $B-$spline bases and knot sequences.}
    \label{}
\end{figure}

The $l$-th $B-$spline function is supported on an interval of the form
$[t_l,t_{l+d+1})$ and is identically zero outside this interval. Consequently,
the index $l$ must satisfy $l+d+1\le m$, and therefore $l=1,\ldots,m-d-1$, yielding exactly the $m-d-1$ $B-$spline functions. Each function corresponds to a distinct group of
$d+2$ consecutive knots. Figure~\ref{fig:preliminaries:Bspline_basis} shows a $B-$spline basis of degree $d = 2$ defined over the equidistant knot sequence $t = \left\{0, 1, 2, 3, 4, 5, 6\right\}$.

For fixed degree $d$ and knot sequence $t$, the associated $B-$spline functions
span a finite-dimensional space of piecewise polynomial functions with that
degree and knot structure. In this mathematical sense, they form a basis of
this space.

Using such a $B-$spline basis of degree $d$, the function $f$ in the regression model defined in~\eqref{eq:preliminaries:univariate_model} can be represented in the domain $[x_1, x_n]$ of observed $X$  values     as a linear combination of the $B-$spline
functions of degree $d$ associated with  a fixed knot sequence. For this, the covariate domain $[x_1, x_n]$ is partitioned into \( k \) intervals, by specifying \( k + 1 \)  knots, which are denoted as internal knots  and may be either equidistant or nonequidistant.  
 Both the degree $d$ and the number of intervals $k$ are parameters chosen by the user. To ensure an accurate estimation and avoid undesirable effects, each interval must be covered by \(d+1\) nonzero $B-$spline functions \citep{wand2008semiparametric}. Thus, using only the internal knots to build the $B-$spline basis can raise issues in the estimation of $f$ due to insufficient coverage at the intervals near the boundaries (see Figure~\ref{fig:preliminaries:Bspline_basis}). To avoid this problem, the knot sequence must be extended beyond the observed data $[x_1,x_n]$ by adding $d$ external knots on each side, ensuring that all internal intervals are covered with enough $B-$spline functions. Therefore, a knot sequence of $m=k+2d+1$ knots is needed to adequately estimate $f$, where $k+1$ of the knots are internal and $2d$ are external. Thus, the extended knot sequence is $\mathbf{t} = \{ t_{q}\}_1^{k+2d+1}$, with $t_{d+1} = x_{1}$ and $ t_{k+d+1} = x_{n}$.  Figure~\ref{fig:preliminaries:Bspline_basis_extended} illustrates the  $B-$spline basis of degree $d=2$ resulting from extending the internal knot sequence  $\mathbf{t}=\{0,1,2,3,4,5,6\}$ with two additional external knots on each side. Compared to the nonextended $B-$spline basis shown in Figure~\ref{fig:preliminaries:Bspline_basis}, this extended version ensures that all internal intervals are fully covered by $d+1=3$ nonzero $B-$spline functions.

Then, function $f$ in~\eqref{eq:preliminaries:univariate_model} can be approximated in the interval $[x_1,x_n]$ with a curve $S$ that is a linear combination of the $k + d$ $B-$spline functions in the $B-$spline basis  of degree $d$ resulting from the extended knot sequence $\mathbf{t} = \{t_q\}_{q=1}^{k+2d+1}$. More precisely, for $x \in [x_1,x_n]$,

\begin{equation}
   S(x)=\sum_{l=1}^{k +d} \theta_{l}B_{l,d,\mathbf{t}}(x),
    \label{eq:preliminaries:f_univariate_estimation}
\end{equation}

\noindent where  $\theta_l$, with $l = 1, \dots, k + d$, are the regression coefficients that need to be estimated, together with the intercept $\alpha$ in~\eqref{eq:preliminaries:univariate_model}. Since $S$ is a linear combination of smooth piecewise polynomial functions, it is itself a smooth piecewise polynomial function.

To estimate $f$ in~\eqref{eq:preliminaries:univariate_model} by $S$ in~\eqref{eq:preliminaries:f_univariate_estimation}, the least squares criterion is commonly used. Then, the regression coefficients $\theta_l$, $l=1,\dots, k+d$, in~\eqref{eq:preliminaries:f_univariate_estimation} and the intercept $\alpha$ in~\eqref{eq:preliminaries:univariate_model} are estimated as the solution of the following quadratic program:

\begin{equation}
    \min_{\alpha, \, \theta_1, \dots, \theta_{k+d}} \sum_{i=1}^{n} \left( y_i - \left(\alpha + \sum_{l=1}^{k+d} \theta_l B_{l,d,\mathbf{t}}(x_i) \right) \right)^2.
    \label{eq:preliminaries:LSE}
\end{equation}

\noindent  The decision variables in the problem~\eqref{eq:preliminaries:LSE} are the regression coefficients $\theta_l$ plus the intercept $\alpha$. Then, the number of decision variables depends only on the number of $B-$spline functions in the $B-$spline basis, determined by the choice of parameters $k$ and $d$, plus one for the intercept. Consequently, the number of decision variables remains independent of the number $n$ of observations. This is a key advantage since it means that the computational complexity of problem~\eqref{eq:preliminaries:LSE} does not scale with the size of the dataset.

Problem~\eqref{eq:preliminaries:LSE} can be written in a more compact way by using the following matrix notation. Let $\mathbf{y} = (y_1, \dots, y_n)^\top$ be the response vector containing the observed values of the response variable. Let $\boldsymbol{\theta}= \left( \alpha, \boldsymbol{\theta}_1^\top\right)^\top = \left(\alpha, \theta_1, \dots, \theta_{k+d}\right)^\top$ be the vector of  coefficients of the model and that need to be estimated. Let 
$\mathbf{B}_1$ be a matrix  of size  $n \times (k+d)$ containing the evaluations of the $k+d$ $B-$spline functions at the observed $X$ values, i.e., $(\mathbf{B}_1)_{il} = B_{l,d,\mathbf{t}}(x_i)$, for $i = 1,\dots,n,\; l = 1,\dots,k+d $. Let $\mathbf{B} = (\mathbf{1} : \mathbf{B}_1)$ be the full design matrix  which organizes the basis evaluations to relate the inputs to the model coefficients. Note that the symbol ``$:$'' denotes the concatenation of matrices by columns. The first column, $\mathbf{1}$, is a vector of ones of length $n$, which provides the baseline response level. Using this matrix notation, problem~\eqref{eq:preliminaries:LSE} becomes

\begin{equation}
    \min_{\boldsymbol{\theta} \in \mathbb{R}^{1+k+d}} \quad \| \mathbf{y} - \mathbf{B} \boldsymbol{\theta} \|^2_2.
    \label{eq:preliminaries:LSE_matrix}
\end{equation}

\noindent Since the basis functions in $\mathbf{B}_1$ sum to one at any point, this raises an identifiability issue when estimating both $\alpha$ and $\boldsymbol{\theta}_1$ simultaneously, and causes the design matrix $\mathbf{B}$ to be rank-deficient, making the least squares problem in~\eqref{eq:preliminaries:LSE_matrix} ill-posed, as the solution is no longer unique. To ensure a unique solution, we  follow the strategy proposed by \citet{wood2020inference} and we address this by enforcing the function $S$ to have zero mean over the observed covariate values, i.e.,  $\sum_{i=1}^n S(x_i) = \mathbf{1}^\top \mathbf{B}_1 \boldsymbol{\theta}_1 = 0$.  To incorporate this into the estimation of the model parameters, we define $\mathbf{P}_1^I = (\mathbf{B}_1^\top \mathbf{1})(\mathbf{1}^\top \mathbf{B}_1)$ and add the quadratic term $\boldsymbol{\theta}_1^\top \mathbf{P}_1^I \boldsymbol{\theta}_1$ to the objective function in~\eqref{eq:preliminaries:LSE_matrix}. This term is mathematically equivalent to $(\mathbf{1}^\top \mathbf{B}_1 \boldsymbol{\theta}_1)^2$. By adding this term to the objective function in~\eqref{eq:preliminaries:LSE_matrix}, the optimization effectively drives the mean of the smooth component toward zero. While this term only affects the spline coefficients $\boldsymbol{\theta}_1$, we formally define the identifiability penalty matrix $\mathbf{P}^I$ as

\begin{equation*}
\mathbf{P}^I =   
\begin{pmatrix}
0 & \mathbf{0}^\top \\
\mathbf{0} & \mathbf{P}_1^I
\end{pmatrix} ,
\label{eq:preliminaries:PI}
\end{equation*}

\noindent  so that the penalty can be written compactly as $\boldsymbol{\theta}^\top P^I \boldsymbol{\theta}$ for the full parameter vector $\boldsymbol{\theta} = (\alpha, \boldsymbol{\theta}_1^\top)^\top$. 
Thus, the problem for estimating $\boldsymbol{\theta}$  becomes the following penalized least squares problem:

\begin{equation}
    \min_{\boldsymbol{\theta} \in \mathbb{R}^{1+k+d}} \quad\| \mathbf{y} - \mathbf{B} \boldsymbol{\theta} \|^2_2 + \boldsymbol{\theta}^\top \mathbf{P}^I \boldsymbol{\theta}.
    \label{eq:preliminaries:LSE_matrix_penalty}
\end{equation}

\noindent The optimal solution of~\eqref{eq:preliminaries:LSE_matrix_penalty} is given by $\widehat{\boldsymbol{\theta}} = \left( \mathbf{B}^\top \mathbf{B} + \mathbf{P}^I \right)^{-1} \mathbf{B}^\top \mathbf{y}$, with $\widehat{\boldsymbol{\theta}} = \left( \hat \alpha, \widehat{\boldsymbol{\theta}_1}^\top\right)^\top.$ Then, $f$ in~\eqref{eq:preliminaries:univariate_model} is estimated by  the curve $S$ in~\eqref{eq:preliminaries:f_univariate_estimation} using the estimated coefficients $\widehat{\boldsymbol{\theta}_1}.$ In particular, letting $\hat f$ denote the estimator of $f$ evaluated at the
observed covariate values, we have
\[
\hat f = \mathbf{B}_1 \widehat{\boldsymbol{\theta}_1}.
\]
 
\noindent As a consequence of $\hat f$ having zero mean due to the identifiability penalty term, the estimated intercept $\hat{\alpha}$
coincides with the sample mean of the response, i.e.,
$\hat{\alpha} = \frac{1}{n} \sum_{i=1}^n y_i$, and $\hat f$ captures only the
deviations from the mean explained by the covariate $X$.

\begin{figure}[h!]
    \centering        \includegraphics[width=0.7\linewidth]{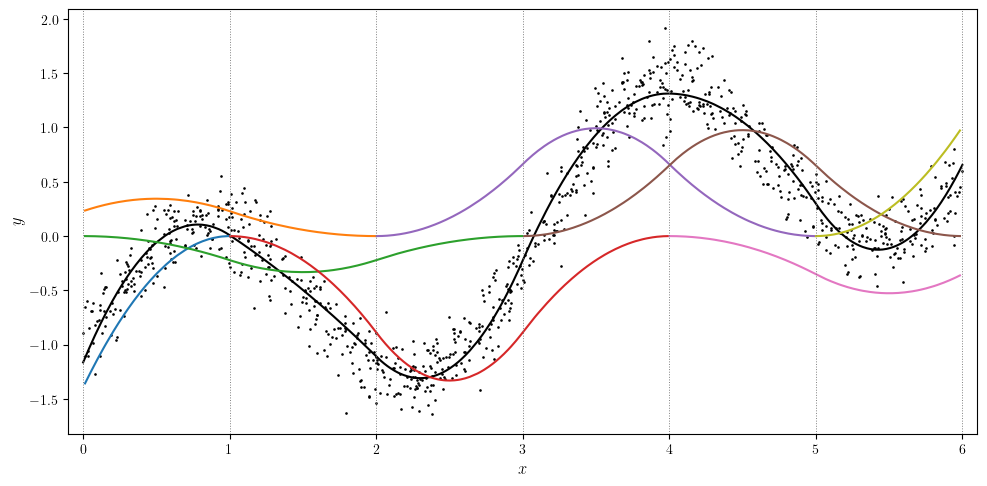}
    \caption{Example of function $f$ estimation (black line) from observed data (dots) using the $B-$spline basis in Figure~\ref{fig:preliminaries:Bspline_basis_extended}. The $B-$spline functions (in colors) are scaled by the regression coefficients, estimated by solving problem~\eqref{eq:preliminaries:LSE_matrix_penalty}.}     \label{fig:preliminaries:f_estimation_B_spline}
\end{figure}

Figure~\ref{fig:preliminaries:f_estimation_B_spline} shows an example of the estimation of the function $f$ in the smooth regression model~\eqref{eq:preliminaries:univariate_model} from a given observed dataset (black dots), using the $B-$spline basis in Figure~\ref{fig:preliminaries:Bspline_basis_extended}. The regression coefficients of  $S$ in~\eqref{eq:preliminaries:f_univariate_estimation} have been estimated by solving problem~\eqref{eq:preliminaries:LSE_matrix_penalty}. The colored curves correspond to the $B$-splines functions, which have been scaled according to these estimated regression coefficients.  Some of them appear inverted as a result of negative coefficient values.

\subsection{Shape-constrained smooth univariate regression model}
\label{sec:preliminaries:univariate_shape_constraints}

In many real-world applications, there is often a need to estimate regression models that not only fit the data well, but also follow certain expected patterns, such as in terms of their shape or sign. Shape-constrained regression helps with this by enforcing these properties during estimation, ensuring that the result aligns with theoretical or expert knowledge expectations. Some authors address this challenge with sum of squares estimators \citep{curmei2025shape}. In this work, we adapt the approach developed by \cite{Navarro2023AMC} to impose shape constraints in the function $S$ given by~\eqref{eq:preliminaries:f_univariate_estimation}
that approximates $f$ in the regression model~\eqref{eq:preliminaries:univariate_model} as described in Section~\ref{sec:preliminaries:univariate_regression}.

In the following, we explain how to impose bound constraints, i.e., restricting the curve $S$ to lie within a specified lower and upper bound across its entire domain $[x_1,x_n]$.  Notice that monotonicity or curvature shape constraints can also be easily imposed by requiring the first or second derivatives of $S$ to satisfy sign constraints. Since $S$  is modeled as a smooth piecewise polynomial function, its derivatives are readily available, making it easy to impose such constraints. In addition, these different shape requirements can be combined by adding the corresponding set of constraints to the resulting problem.

Suppose that it is known from expert knowledge that the response variable $Y$ in the regression model in~\eqref{eq:preliminaries:univariate_model}  lies within a fixed interval  $[L,U]$, for some $L,U \in \mathbb{R}$ with $L<U$, in the observed domain $[x_1,x_n]$ of the  covariate $X$. Therefore, to remain true to this prior knowledge,  the  regression model must also satisfy these bounds. That is,

\begin{equation}
L \leq  \alpha +  f(x) \leq U, \quad \forall x \in [x_1, x_n],
\label{eq:preliminaries:univariate_bounded_model}
\end{equation}

\noindent where  $\alpha$ and $f$ are estimated as in Section~\ref{sec:preliminaries:univariate_regression}. Since the estimation of $\alpha$ is a known constant equal to the mean of the observed responses, the conditions in~\eqref{eq:preliminaries:univariate_bounded_model} translate into bounds on $ f$, namely,

\begin{equation}
L -  \alpha \leq   f(x) \leq U - \alpha, \quad \forall x \in [x_1, x_n]. 
\label{eq:preliminaries:f_bounded}
\end{equation}

\noindent Therefore, to estimate $f$ under the bounds in~\eqref{eq:preliminaries:f_bounded},  $S$ in~\eqref{eq:preliminaries:f_univariate_estimation} must satisfy the same bounds. Note that these bounds should only be imposed when there is prior expert knowledge that the response variable lies in the interval $[L,U]$. Otherwise, inappropriate bounds may lead to an erroneous approximation. When such bound information is available, an upper bound of $U-L$ on the approximation error of the fitted regression model over the entire covariate domain is obtained. However, this bound is normally loose and should only be considered as a naive worst-case bound and not as an informative approximation guarantee.

To estimate $ f$ under the bound constraints in~\eqref{eq:preliminaries:f_bounded} through the curve $S$ in~\eqref{eq:preliminaries:f_univariate_estimation}, we adapt the approach in \cite{Navarro2023AMC} for penalized splines ($P-$splines) to the case of univariate shape-constrained smooth regression with $B-$splines. Their method enforces shape constraints in $S$ by leveraging the necessary and sufficient conditions for nonnegativity of an univariate polynomial over an interval, as stated in \cite{bertsimas2002relation}.  In particular, polynomial univariate nonnegativity over an interval is equivalent to a sum-of-squares representation and \cite{bertsimas2002relation} provide necessary and sufficient conditions for this representation to hold.  Since $S$ is a smooth piecewise polynomial, these conditions are applied locally to each of the $k$ internal intervals within the data domain $[x_1,x_n].$

Then, the lower and upper bounds on $ f$ as stated in~\eqref{eq:preliminaries:f_bounded} are enforced by requiring nonnegativity over the domain $[x_1,x_n]$ of the auxiliary functions $f_L(x) \coloneq  f(x) - (L -   \alpha)$ and $f_U(x)\coloneq -f(x)  + ( U -  \alpha)$.  As detailed in  \cite{Navarro2023AMC}, enforcing the aforementioned necessary and sufficient polynomial nonnegativity conditions in each case and to each internal interval within the data domain $[x_1,x_n]$ leads to new decision variables and constraints that are incorporated into problem~\eqref{eq:preliminaries:LSE_matrix_penalty} to estimate $\boldsymbol{\theta}_1$ in $S$. The following conic optimization problem is obtained, which corresponds to the formulation proposed in~\cite{Navarro2023AMC}, adapted from $P$-splines to the $B$-spline representation considered in this work:

\begin{equation*}
\begin{aligned}
    &  \min   \quad \| \mathbf{y} - \mathbf{B} \boldsymbol{\theta} \|^2_2 + \boldsymbol{\theta}^\top \mathbf{P}^I \boldsymbol{\theta}  \\
    & \text{ s.t.}\\
     &  \quad   \langle \mathbf{H}_{\ell}, \mathbf{Z}^L_q \rangle_F = 0,    \ \quad \quad  \quad  \quad      q = d+1,\dots,d+k, \  \ell = 1,\dots,d, \\
    &  \quad  \langle \mathbf{H}_{\ell}, \mathbf{Z}^U_q \rangle_F = 0, \   \quad   \quad  \quad  \quad     q = d+1,\dots,d+k, \   \ell = 1,\dots,d,\\
    &  \quad   \begin{pmatrix}
        \langle \mathbf{H}_{d+1}, \mathbf{Z}^L_q \rangle_F \\
        \vdots\\
        \langle \mathbf{H}_{2d+1}, \mathbf{Z}^L_q \rangle_F
    \end{pmatrix}
    =  \mathbf{W}_q \mathbf{G}_q \boldsymbol{\theta}_1- \mathbf{W}_q
    \begin{pmatrix}
        L -  \hat \alpha \\
        0\\
        \vdots\\
        0
    \end{pmatrix}, 
     \ \quad     q  = d+1, \dots, d+k,  \\
    &  \quad   \begin{pmatrix}
        \langle \mathbf{H}_{d+1}, \mathbf{Z}^U_q \rangle_F \\
        \vdots\\
        \langle \mathbf{H}_{2d+1}, \mathbf{Z}^U_q \rangle_F
    \end{pmatrix}
    =  -\mathbf{W}_q \mathbf{G}_q \boldsymbol{\theta}_1+ \mathbf{W}_q
    \begin{pmatrix}
        U -  \hat \alpha \\
        0\\
        \vdots\\
        0
    \end{pmatrix}, 
     \ \quad     q  = d+1, \dots, d+k,\\
    &  \quad  \mathbf{Z}^U_q, \mathbf{Z}^L_q \in \mathcal{S}_+^{d+1}, \quad  \quad  \quad  \quad    q  = d+1, \dots, d+k,\\
    &  \quad  \boldsymbol{\theta} = \left( \alpha,  \boldsymbol{\theta}_1^\top\right)^\top \in \mathbb{R}^{1+k+d} ,
\end{aligned}
\label{eq:preliminaries:univariate_shape_constrained_problem}
\end{equation*}

\noindent where $\langle \cdot, \cdot \rangle_F$ denotes the Frobenius inner product. The constraints involving superscript $U$ enforce the upper bound on $f$, while those with superscript $L$ enforce the lower bound. The matrices  $\mathbf{Z}^U_q$ and $\mathbf{Z}^L_q$, with  $q  = d+1, \dots, d+k$, are decision variables that arise when applying the necessary and sufficient conditions for polynomial nonnegativity 
to the approximation of the auxiliary functions $ f_L(x)$  and $ f_U(x)$ at each internal interval $[t_q,t_{q+1})$ within the domain $[x_1, x_n]$.  The set $\mathcal{S}_+^{d+1}$ denotes the cone of positive semidefinite matrices of order $d+1$. Furthermore, the matrices  $\mathbf{H}_\ell, \ \ell=1,\dots, 2d+1$,  $\mathbf{W}_q, \ q=d+1,\dots,d+k$ and  $\mathbf{G}_q, \ q=d+1,\dots,d+k$, also arise from the application of the nonnegativity conditions and help in representing the set of linear constraints compactly. They are described next. Matrices $\mathbf{H}_\ell \in \mathbb{R}^{(d+1) \times (d+1)}$, for $\ell = 1, \dots, 2d+1$, have entries defined as

\begin{equation*}
        (\mathbf{H}_\ell)_{ij} = 
\begin{cases}
1 & \text{if } i + j = 2(d+1-l)+1,\quad \text{for } l = 1,\dots,d, \\
1 & \text{if } i + j = 2(2d+2-l), \quad \text{for } l =  d+1,\dots, 2d+1 , \\
0 & \text{otherwise}.
\end{cases}
\label{eq:preliminaries:H_matrices}
\end{equation*}

\noindent In addition, for each internal interval $[t_q,t_{q+1})$ within the data domain $[x_1, x_n]$, with $q = d+1, \dots, d+k$, a matrix  $\mathbf{W}_q \in \mathbb{R}^{(d+1) \times (d+1)}$ is defined as 

\begin{equation*}
 (\mathbf{W}_q)_{ij} = \sum_{m=\max(0,i+j-2-d)
}^{\min(i-1,j-1)} \binom{j-1}{m}\binom{d-j+1}{i-1-m} t_q^{j-1-m} t_{q+1}^{m}.
\label{eq:preliminaries:W_matrices}
\end{equation*}

\noindent Similarly, for each internal interval $[t_q,t_{q+1})$, with $q = d+1, \dots, d+k$, a matrix  $\mathbf{G}_q \in \mathbb{R}^{(d+1) \times (d+k)}$ is defined as 

\begin{equation*}
\mathbf{G}_q = \left( 
\mathbf{0}_{(d+1) \times (q - (d+1))} \  {:} \ 
\mathbf{{\tilde{G}}}_q \ {:} \ 
\mathbf{0}_{(d+1) \times (k + d - q)} 
\right),
\label{eq:preliminaries:G_matrices}
\end{equation*}

\noindent where $\mathbf{\tilde{G}}_q \in \mathbb{R}^{(d+1)\times (d+1)}$ is a matrix containing the polynomial coefficients of each of the $d+1$ pieces of $B-$spline functions that overlap in the interval $[t_q,t_{q+1})$ (see Figure~\ref{fig:preliminaries:Bspline_basis_extended}). More specifically, $  ({\mathbf{\tilde{G}}}_q)_{ij} = g_{q-d+j-1,q,i}$,  being $g_{q-d+j-1,q,i}$  the polynomial coefficient of degree $i$ (with $i=0,1,\dots,d$) of the segment of the $B-$spline function starting at knot  $t_{q-d+j-1}$ restricted to the interval  $[t_q, t_{q+1})$. 
All these previous concepts are necessary for Section~\ref{sec:shape_constrained_regression}.  Their technical details, as well as the derivation and justification of the optimization problem above, are further explained in \cite{Navarro2023AMC}.

\subsection{Smooth additive regression models}
\label{sec:preliminaries:smooth_additive_regresion}

This section extends the $B-$spline estimation approach for unconstrained smooth univariate regression models, described in Section~\ref{sec:preliminaries:univariate_regression}, to smooth additive regression models. In these models, multiple covariates influence the response variable in a smooth and additive manner, offering flexible modeling while maintaining interpretability \citep{hastie2017generalized}.

Formally, consider a set of $n$ observations $\{(\mathbf{x}_i, y_i), i = 1, \dots, n \}$ drawn from $p$ continuous covariates $X_1, \dots, X_p$, i.e., $\mathbf{x}_i=(x_{i1}, x_{i2}, \dots,  x_{ip})$ and a continuous response variable $Y$, respectively. The smooth additive regression model is defined as

    \begin{equation}
        y_i = \alpha + f_1(x_{i1}) + \dots +  f_p(x_{ip}) + \epsilon_i, \quad i = 1, \dots, n,    \label{eq:preliminaries:smooth_additive_model}
    \end{equation}

\noindent where $\epsilon_i \in \mathbb{R}$ are independent error terms with zero mean ($i =1,\dots,n$),  $\alpha \in \mathbb{R}$ and each $f_j:~ [\underline{x}_{j}, \overline{x}_{j}] \subset \mathbb{R} \rightarrow \mathbb{R}$, for $j=1,\dots, p$,  is a univariate smooth function defined in the range of observed values of $X_j$. Here, $\underline{x}_{j}$ and $\overline{x}_{j}$ denote the minimum and maximum observed values of $X_j$, respectively.

The goal is to estimate the intercept $\alpha$ along with the smooth functions $f_j$, $j=~1,\dots,p$,  in~\eqref{eq:preliminaries:smooth_additive_model}  using the observed covariate values to closely approximate  the response values. Following the methodology described in Section~\ref{sec:preliminaries:univariate_regression}, we adopt a $B-$spline approach. For this, the domain of each variable $X_j$ is divided into a fixed $k_j$ number of intervals by specifying $k_j +1$ internal knots, which may be either equidistant or not. The degree of the $B-$spline basis for each variable $X_j$  is also fixed at $d_j$, for $j=1, \dots, p$. This setup yields an extended knot sequence for each covariate $X_j$, namely, $\mathbf{t}_j = \{t_{jq}\}_{q=1}^{k_j+2d_j+1}$, with $t_{j(d_j+1)}=\underline{x}_j$ and $t_{j(k_j+d_j+1)}=\overline{x}_j$. Consequently, there are $k_j + d_j$ $B-$spline functions to estimate each smooth function $f_j$ associated with variable $X_j$, $j=1,\dots, p$.  Each function $f_j$ is then approximated by a curve $S_j$ that is a
a linear combination of them. More precisely, for $x \in [\underline{x}_j,\overline{x}_j]$, $S_j(x)$ is given by

\begin{equation}
    S_j(x)=\sum_{l=1}^{k_j + d_j} \theta_{jl}B_{l,d_j,\mathbf{t}_j}(x), 
    \label{eq:preliminaries:fj_estimation}
\end{equation}

\noindent where $\theta_{jl}$, with $l = 1,\dots, k_j + d_j$, are the regression coefficients associated with the $B-$spline functions  of the variable $X_j$. 

The estimation of each $f_j$ in  the smooth additive model in~\eqref{eq:preliminaries:smooth_additive_model} requires estimating the regression coefficients of its corresponding  curves $S_j, \ j=1,\dots,p$. Thus, the model in~\eqref{eq:preliminaries:smooth_additive_model} is estimated by solving the following penalized least squares problem, which is similar to~\eqref{eq:preliminaries:LSE_matrix_penalty},
but extended to the multivariate additive case:

\begin{equation}
    \min_{\boldsymbol{\theta} \in \mathbb{R}^{1+\sum_{j=1}^p k_j+d_j}} \quad\| \mathbf{y} - \mathbf{B} \boldsymbol{\theta} \|^2_2 + \boldsymbol{\theta}^\top \mathbf{P}^I \boldsymbol{\theta},
    \label{eq:preliminaries:PLSE_additive_matrix}
\end{equation}

\noindent where $\mathbf{y} = (y_1, \dots, y_n)^\top$ is the response vector. The vector $\boldsymbol{\theta}$ and the matrices $\mathbf{B}$ and  $\mathbf{P}^I$ are extended to incorporate the $p$ covariates. 
The full parameter vector to be estimated is $\boldsymbol{\theta} = (\alpha, \boldsymbol{\theta}_1^\top, \dots, \boldsymbol{\theta}_p^\top)^\top$, where each $\boldsymbol{\theta}_j=(\theta_{j1}, \dots, \theta_{j(k_j + d_j)})^\top$ contains the regression coefficients for the $B-$spline functions corresponding to covariate $X_j$, as explained in ~\eqref{eq:preliminaries:fj_estimation}. The full design matrix is now given by $\mathbf{B} =  (\mathbf{1} : \mathbf{B}_1: \dots : \mathbf{B}_p )$,
where $\mathbf{1}$ is a vector of ones of length $n$ and $\mathbf{B}_j$ is an $n \times (k_j + d_j)$ matrix associated with covariate $X_j$ containing the evaluations of its $k_j + d_j$ $B-$spline functions at the observed  values. More specifically, $(\mathbf{B}_j)_{il} = B_{l,d_j,\mathbf{t}_j}(x_{ij})$, for $i = 1,\dots,n$ and $l = 1,\dots,k_k+d_j$. The identifiability penalty matrix $\mathbf{P}^I$ in the multivariate setting  is a  block diagonal matrix defined as

\begin{equation*}
\mathbf{P}^I =
\begin{pmatrix}
0 & \mathbf{0}^\top & \cdots &  \mathbf{0}^\top\\
\mathbf{0} & \mathbf{P}_1^I & \cdots & \mathbf{0}^\top\\
\vdots & \vdots& \ddots & \vdots\\
\mathbf{0} & \mathbf{0} & \cdots & \mathbf{P}_p^I 
\end{pmatrix},
\label{eq:preliminaries:PI_additive_model}
\end{equation*}

\noindent where $\mathbf{P}_j^I = \mathbf{B}_j^\top \mathbf{1} \mathbf{1}^\top \mathbf{B}_j$, for $j = 1, \dots, p$, is the identifiability penalty matrix of covariate $X_j$ and $\mathbf{0}$ denotes a zero vector of appropriate length for dimension  compatibility. Recall that  the explicit optimal solution of problem~\eqref{eq:preliminaries:PLSE_additive_matrix} is $\widehat{\boldsymbol{\theta}} =~ \left( \mathbf{B}^\top \mathbf{B} + \mathbf{P}^I \right)^{-1} \mathbf{B}^\top \mathbf{y}$, where in this case $\widehat{\boldsymbol{\theta}} = \left( \hat \alpha, \widehat{\boldsymbol{\theta}_1}^\top,\dots, \widehat{\boldsymbol{\theta}_p}^\top\right)^\top.$ Then, each $f_j$ is estimated by $S_j$ in~\eqref{eq:preliminaries:fj_estimation}  with its estimated parameter $\widehat{\boldsymbol{\theta}_j}$, i.e., $\hat f_j = \mathbf{B}_j \widehat{\boldsymbol{\theta}_j}$. 
Recall as well that the identifiability penalty forces each estimated function $\hat{f}_j$ to have zero mean over the observed values. As a consequence, the estimated intercept $\hat{\alpha}$ is the mean of the observed responses.

\section{Smooth additive regression model with constraints}
\label{sec:shape_constrained_regression}

The estimation of smooth additive regression models, as described in Section~\ref{sec:preliminaries:smooth_additive_regresion}, does not ensure that the fitted model satisfies any shape requirement. However, such requisites are often pivotal for reflecting theoretical or  expert knowledge. This limitation motivates the development of a novel shape-constrained estimation methodology for smooth additive models, which also allows the incorporation of other pointwise constraints easily.

Consider the smooth additive regression model in~\eqref{eq:preliminaries:smooth_additive_model}. Following the discussion in Section~\ref{sec:preliminaries:univariate_shape_constraints}, suppose that  the response variable $Y$ lies within an interval $[L,U]$ over the covariates domain. Then, the regression model must satisfy

\begin{equation}
L \leq  \alpha +  f_1(x_1) +\cdots +  f_p(x_p) \leq U, \quad \forall (x_1, \dots, x_p) \in [\underline{x}_1,\overline{x}_1]\times \dots \times [\underline{x}_p,\overline{x}_p].
\label{eq:shape_constrained_additive:bounds}
\end{equation}

\noindent  Recall that this condition should only be imposed when supported by prior expert knowledge. Since 
$\alpha$ is a constant estimated as the mean of the response values, condition~\eqref{eq:shape_constrained_additive:bounds} is equivalent to

\begin{equation}
L -   \alpha\leq  f_1(x_1) +\cdots +  f_p(x_p) \leq U -   \alpha, \quad \forall (x_1, \dots, x_p) \in [\underline{x}_1,\overline{x}_1]\times \dots \times [\underline{x}_p,\overline{x}_p].
\label{eq:shape_constrained_additive:bounds_sum_fs}
\end{equation}

\noindent Each function $f_j$,  $j=1,\dots, p$, is approximated by a curve $S_j$  given by~\eqref{eq:preliminaries:fj_estimation},  which is a piecewise univariate polynomial. Then the estimation of functions $f_j$ satisfying condition~\eqref{eq:shape_constrained_additive:bounds_sum_fs} involves imposing bounds on a sum of univariate polynomials.  Analogously to the univariate case in~\eqref{eq:preliminaries:f_bounded}, enforcing~\eqref{eq:shape_constrained_additive:bounds_sum_fs} is equivalent to requiring nonnegativity of the auxiliary expressions
$\sum_{j=1}^p f_j(x_j) - (L-\alpha)$ and
$-\sum_{j=1}^p f_j(x_j) + (U-\alpha)$
over the Cartesian product domain of the variables. However, to the authors’ knowledge, there is no direct extension of the necessary and sufficient conditions in \cite{bertsimas2002relation} for the nonnegativity of univariate polynomials over an interval to the additive setting involving sums of polynomials (see \cite{blekherman2012semidefinite}). Therefore, satisfying condition~\eqref{eq:shape_constrained_additive:bounds_sum_fs} is not straightforward. 

In this work, we address condition~\eqref{eq:shape_constrained_additive:bounds_sum_fs} by bounding each individual function $ f_j$ as follows:

\begin{equation}
   \omega^L_j (L- \alpha) \leq   f_j(x) \leq \omega^U_j (U- \alpha), \quad  j = 1,\dots,p, \quad \forall x \in [\underline{x}_j,\overline{x}_j],
\label{eq:shape_constrained_additive:condition_fj}
\end{equation}

\noindent where $\omega^L_j,\omega^U_j \in (0,1), \  j=1,\dots,p$, $\sum_{j=1}^p \omega^L_j =1$ and $\sum_{j=1}^p \omega^U_j =1$. 
The componentwise bounds in~\eqref{eq:shape_constrained_additive:condition_fj} guarantee that the sum of the functions satisfies the global bound in~\eqref{eq:shape_constrained_additive:bounds_sum_fs}: summing~\eqref{eq:shape_constrained_additive:condition_fj} over each function $f_j$ yields
\[
(L-\alpha)\sum_{j=1}^p \omega_j^L \;\le\; \sum_{j=1}^p f_j(x_j) \;\le\; (U-\alpha)\sum_{j=1}^p \omega_j^U,
\]
\noindent and since $\sum_{j=1}^p\omega_j^L=\sum_{j=1}^p\omega_j^U=1$, ~\eqref{eq:shape_constrained_additive:bounds_sum_fs} is obtained. Hence, the conditions in~\eqref{eq:shape_constrained_additive:condition_fj} are sufficient (but not necessary) to enforce~\eqref{eq:shape_constrained_additive:bounds_sum_fs}. Notice that the approach in ~\eqref{eq:shape_constrained_additive:condition_fj}  is possible because of the identifiability constraint that forces the estimation of each function $ f_j$ to be centered at zero (see Section~\ref{sec:preliminaries:smooth_additive_regresion}). This centering plays a crucial role  since it makes the componentwise bounds in~\eqref{eq:shape_constrained_additive:condition_fj} natural. Each $f_j$ is allowed to vary around zero as much as necessary to adapt to the data while satisfying the bounds in~\eqref{eq:shape_constrained_additive:condition_fj}.

To estimate the weights $\omega^L_j$ and $\omega^U_j$, $j=1,\dots,p$, we propose a general data-driven approach. For this, an unconstrained smooth additive regression model as in~\eqref{eq:preliminaries:PLSE_additive_matrix} is first fitted.   Then, we obtain the fitted values of each univariate function $f_j$ in the model. These are the evaluations of the fitted univariate functions  $\widehat f_j$  at each observed covariate value $x_{ij}$, i.e., $\widehat f_j(x_{ij})$ for $i=1,\ldots,n$ and $j=1,\ldots,p$. We then compute its minimum and maximum values: $\underline{\hat f_j}=~\min(\widehat f_j(x_{ij}))_{i=1}^n$ and $ \overline{\hat f_j}=~\max(\widehat f_j(x_{ij}))_{i=1}^n,$  $j=1,\dots,p$. The lower and upper weights, namely $\omega_j^L$ and $\omega_j^U,$ respectively, are computed by normalizing these minima and maxima across all covariates, yielding

\begin{equation*}
        \omega^L_j  = \dfrac{ \underline{\hat f_j}}{\sum_{j=1}^p  \underline{\hat f_j}}, \quad         \omega^U_j  = \dfrac{ \overline{\hat f_j}}{\sum_{j=1}^p  \overline{\hat f_j}},  \quad j=1,\dots,p.
\end{equation*}

\noindent This choice of the weights  $\omega^L_j$ and $\omega^U_j$, $j=1,\dots,p$, using the minima and maxima of the fitted values $\widehat f_j(x_{ij})$, $i=1,\ldots,n$ calibrates the componentwise bounds in~\eqref{eq:shape_constrained_additive:condition_fj} accounting for differences in the ranges of the fitted values of each $f_j$. If expert knowledge about the influence of each covariate on the response is available, the weights $\omega^L_j$ and $\omega^U_j$, $j=1,\dots,p$ may be determined based on that knowledge. 

 Building on the discussion above, to estimate a smooth additive regression model satisfying the global bound constraint~\eqref{eq:shape_constrained_additive:bounds_sum_fs}, we enforce the componentwise bounds~\eqref{eq:shape_constrained_additive:condition_fj} for each univariate function $f_j$. Since these are univariate constraints, the bound-constrained estimation of each $f_j$ can be carried out using the shape-constrained univariate regression methodology of \cite{Navarro2023AMC}, detailed in Section~\ref{sec:preliminaries:univariate_shape_constraints}. Hence, a set of constraints for the lower and upper bounds of each $f_j$ is added to the smooth additive regression model estimation~\eqref{eq:preliminaries:PLSE_additive_matrix}, leading to the following formulation:

{ \fontsize{8.5pt}{11pt}\selectfont
\begin{equation}
\begin{aligned}
    & \min  \quad   \| \mathbf{y} - \mathbf{B} \boldsymbol{\theta} \|^2_2 + \boldsymbol{\theta}^\top \mathbf{P}^I \boldsymbol{\theta}  \\
    & \text{ s.t.}\\
    &   \quad \langle \mathbf{H}_{j\ell_j}, \mathbf{Z}^L_{jq_j} \rangle_F = 0, \ \quad \quad \quad  j  = 1, \dots, p, \   q_j  = d_j+1, \dots, d_j+k_j, \  \ell_j = 1, \dots, d_j, \\
    & \quad  \langle \mathbf{H}_{j\ell_j}, \mathbf{Z}^U_{jq_j} \rangle_F = 0, \ \quad \quad \quad   j  = 1, \dots, p, \   q_j  = d_j+1, \dots, d_j+k_j, \  \ell_j = 1, \dots, d_j, \\
    & \quad    \begin{pmatrix}
        \langle \mathbf{H}_{j(d_j+1)}, \mathbf{Z}^L_{jq_j} \rangle_F \\
        \vdots\\
        \langle \mathbf{H}_{j(2d_j+1)}, \mathbf{Z}^L_{jq_j} \rangle_F
    \end{pmatrix}
    =  \mathbf{W}_{jq_j} \mathbf{G}_{jq_j} \boldsymbol{\theta}_j- \mathbf{W}_{jq_j}
    \begin{pmatrix}
       \omega^L_j (L -  \hat \alpha) \\
        0\\
        \vdots\\
        0
    \end{pmatrix}, 
    \ \begin{array}{l}
      j = 1, \dots, p \\
     q_j = d_j+1, \dots, d_j+k_j
\end{array}, \\
    & \quad   \begin{pmatrix}
        \langle \mathbf{H}_{j(d_j+1)}, \mathbf{Z}^U_{jq_j} \rangle_F \\
        \vdots\\
        \langle \mathbf{H}_{j(2d_j+1)}, \mathbf{Z}^U_{jq_j} \rangle_F
    \end{pmatrix}
    =  -\mathbf{W}_{jq_j} \mathbf{G}_{jq_j} \boldsymbol{\theta}_j+ \mathbf{W}_{jq_j}
    \begin{pmatrix}
         \omega^U_j(U -  \hat \alpha) \\
        0\\
        \vdots\\
        0
    \end{pmatrix}, 
    \ \begin{array}{l}
     j = 1, \dots, p \\
     q_j = d_j+1, \dots, d_j+k_j
\end{array}, \\
    &\quad  \mathbf{Z}^U_{jq_j}, \  \mathbf{Z}^L_{jq_j} \in \mathcal{S}_+^{d_j+1}, \ \quad \quad \quad  j=1,\dots,p, \  q_j  = d_j+1, \dots, d_j+k_j,\\
    &\quad  \boldsymbol{\theta} = \left( \alpha,  \boldsymbol{\theta}_1^\top,\dots,\boldsymbol{\theta}_p^\top\right)^\top \in \mathbb{R}^{1+\sum_{j=1}^pk_j+d_j},
\end{aligned}
\label{eq:shape_constrained_additive:shape_constrained_problem}
\end{equation}
}

\noindent where $\mathbf{y}$, $\mathbf{B}$, $\mathbf{P}$ and $\boldsymbol{\theta}$ are those introduced in Section~\ref{sec:preliminaries:smooth_additive_regresion}. Matrices $\mathbf{Z}^U_{jq_j}$,  $\mathbf{Z}^L_{jq_j}$, $\mathbf{W}_{jq_j}$, $\mathbf{G}_{jq_j}$ and $\mathbf{H}_{j\ell_j}$, with $j=1,\dots,p$,  $q_j  = d_j+1, \dots, d_j+k_j$ and $\ell_j  = 1, \dots, 2d_j+1$, are defined as in in Section~\ref{sec:preliminaries:univariate_shape_constraints},  adapted for the degree $d_j$ and the number of internal intervals $k_j$ of each covariate.

The proposed approach provides a practical decomposition of the global bound requirement~\eqref{eq:shape_constrained_additive:bounds_sum_fs} into the componentwise constraints~\eqref{eq:shape_constrained_additive:condition_fj}. Therefore, any feasible solution of~\eqref{eq:shape_constrained_additive:shape_constrained_problem} necessarily satisfies~\eqref{eq:shape_constrained_additive:condition_fj} for all function $f_j$ and, consequently, satisfies the global bound~\eqref{eq:shape_constrained_additive:bounds_sum_fs}.

For simplicity, we have focused on explaining the details for model estimation with bound constraints,  but the proposed methodology is easily extended to monotonicity and curvature. These can be imposed by requiring the first or second derivatives of each $S_j$, given by~\eqref{eq:preliminaries:fj_estimation} and approximating $f_j$, to be nonnegative across their respective domains. Since each $S_j$ estimated this way 
will satisfy the desired shape property, their sum will also inherit it. Different shape requirements can be also combined by adding the corresponding constraints to the problem.

In addition to the shape and bound constraints described above, the proposed methodology also allows for the straightforward inclusion of pointwise constraints on the estimated model. These include interpolation constraints (forcing the model to take specific values at specific points), as well as underestimation or overestimation (forcing the model to lie below or above certain values at specific points). Formally, suppose that, at a given subset $\mathcal{X}_A$ of observed covariate values,  $\mathcal{X}_A = \{(x_{i1},x_{i2},\dots, x_{ip}), \ i \in A \subseteq  \{1,\dots,n\} \}$, the estimated model must take values equal to, lower, or higher than the corresponding observed $\mathbf{y}_A$ values. An interpolation, underestimation, or overestimation constraints are readily incorporated into the optimization problem~\eqref{eq:shape_constrained_additive:shape_constrained_problem} by adding the corresponding linear constraints among

\begin{equation*}
\begin{aligned}
        &\mathbf{B}_A\boldsymbol{\theta}  =\mathbf{y}_A,\\
        &\mathbf{B}_A\boldsymbol{\theta}  \leq\mathbf{y}_A,\\
        &\mathbf{B}_A\boldsymbol{\theta}  \geq\mathbf{y}_A,
\end{aligned}
\label{eq:shape_constrained_additive:interpolation_constraints}
\end{equation*}

\noindent where $ \mathbf{B}_A$ is the full design matrix $ \mathbf{B}$ restricted to $\mathcal{X}_A$.

\section{Mixed-Integer Smoothing Surrogate Optimization with Constraints (MISSOC)}
\label{sec:surrog}

In this section, we introduce our method to formulate and solve surrogate MINLPs, the Mixed-Integer Smoothing Surrogate Optimization with Constraints (MISSOC) 
algorithm. First, in Section~\ref{sec:surrog_surrogate-problem} we introduce a mathematical optimization problem as an abstract MINLP and its surrogate counterpart based on the approach developed in Section~\ref{sec:shape_constrained_regression}. Then,  in Section~\ref{sec:surrog_formulation}, we develop a formulation for this surrogate,  that is inspired by the classic multiple choice formulation \citep{vielma2009mixed} and yields an MINLP. In Section~\ref{sec:surrog_the-overall-algorithm}, we present the overall methodology for obtaining a solution to the original problem by solving the the resulting surrogate MINLP. 

\subsection{The surrogate problem}
\label{sec:surrog_surrogate-problem}

Assume that we are given a  mathematical optimization problem as an MINLP, say~\eqref{sec:surr_problem:P}, of the following form:

\begin{equation}
    \begin{aligned}
& \min \quad g_0(\mathbf{x})  \\ 
 & \text{ s.t.}\\
& \quad g_m(\mathbf{x}) \leq 0, \quad \quad  \quad m \in \{ 1,\dots,M\}, \\
& \quad x_j \in \mathbb{Z}, \quad \quad \quad \quad  j  \in I \subseteq \{1,\dots,p\},\\
& \quad \underline{x}_j \leq x_j \leq \overline{x}_j, \quad \quad  j  \in \{1,\dots,p\},\\
& \quad \mathbf{x} \in \mathbb{R}^p,    \end{aligned}
\tag{$P$}
\label{sec:surr_problem:P}
\end{equation}

\noindent where $g_m:\mathbb{R}^p \to \mathbb{R}$, $m=0,\ldots,M$, are scalar functions,  with $g_0(\mathbf{x})$ representing the objective function and $g_m$, $m=1,\ldots,M,$ the constraints. 
$M$ is the number of constraints, $I$ is the index set of integer variables, and $\underline{x}_j$ and $\overline{x}_j$ are the lower and upper bound on variable $j$, respectively. Note that $\underline{x}_j$ and $\overline{x}_j$ can take values $-\infty$ and $+\infty$, respectively, when $x_j$ is not bounded below/above. However, the lower and upper bounds are assumed to take a finite value for each variable appearing nonlinearly in the objective function or constraints. Also note that in MINLPs one usually only assumes that functions $g_m(\mathbf{x})$, $m=0,\ldots,M$, are twice-continuously differentiable and no particular assumption on their convexity is made. Given~\eqref{sec:surr_problem:P}, a surrogate problem can be built as stated in Definition~\ref{def:surrogate}.

\begin{definition}\label{def:surrogate}
    Let $C~\subset~\{0,\ldots,M\}$  be the index set of complicating functions in $(P)$ . Then, the surrogate problem associated with~\eqref{sec:surr_problem:P}, denoted by \eqref{sec:surr_problem:tildeP}, is defined as 
    \begin{equation}
    \begin{aligned}
& \min \quad {\tilde g}_0(\mathbf{x})  \\ 
& \text{ s.t.}\\
& \quad {\tilde g}_m(\mathbf{x}) \leq 0,  \quad \quad  \quad  m \in \{ 1,\dots,M\}, \\
& \quad x_j \in \mathbb{Z},  \quad \quad \quad  \quad  j  \in I \subseteq \{1,\dots,p\},\\
& \quad \underline{x}_j \leq x_j \leq \overline{x}_j, \quad \quad   j  \in \{1,\dots,p\},\\
& \quad \mathbf{x} \in \mathbb{R}^p,    \end{aligned}
\tag{$\tilde{P}$}
\label{sec:surr_problem:tildeP}
\end{equation}
being $\tilde g_m(\mathbf{x}) = g_m(\mathbf{x})$ for $m \in \{0,\dots,M\} \setminus C$, and $\tilde g_m(\mathbf{x})$ the approximation of $g_m(\mathbf{x})$ for $m\in C.$
\end{definition}

\noindent Let us now discuss the role of~\eqref{sec:surr_problem:tildeP}. The main idea is to replace functions in~\eqref{sec:surr_problem:P} that are difficult to handle computationally by approximating functions that are easier. In this sense,  a surrogate of~\eqref{sec:surr_problem:P}  could be obtained by substituting only the subset of the complicating functions and the non complicating functions can remain as they are.  Depending on the properties of the  $\tilde g_m(\mathbf{x})$ functions,  $m\in C$, \eqref{sec:surr_problem:tildeP} could be: i) a relaxation  of~\eqref{sec:surr_problem:P}, when $\tilde g_m(\mathbf{x}) \leq g_m(\mathbf{x})$; ii) a restriction  of~\eqref{sec:surr_problem:P}, when $\tilde g_m(\mathbf{x}) \geq g_m(\mathbf{x})$; iii) an approximation  of~\eqref{sec:surr_problem:P}, otherwise.

Many existing MINLP algorithms rely on solving  a problem~\eqref{sec:surr_problem:P} via a surrogate like \eqref{sec:surr_problem:tildeP}. These approaches can be classified in exact and heuristic methods. The former, which guarantee to provide a global optimum, are mainly based on solving relaxations obtained by convexifying/linearizing functions $g_m(\mathbf{x}),$ ${m\in C}$. Classical examples include spatial branch-and-bound and its variants, see, for example, \cite{belotti2009branching,ryoo-sahinidis-1996,sahinidis1996baron,Smith-Pantelides-1997,Tawarmalani-Sahinidis-2002,Tawarmalani-Sahinidis-2004}. Their drawbacks are that they can be time-consuming and cannot scale up to large instances. In contrast, heuristic algorithms aim to find a good-quality solution in a short amount of time. They are often based on solving surrogate problems~\eqref{sec:surr_problem:tildeP} which are restrictions or approximations of~\eqref{sec:surr_problem:P}. In particular, piecewise linear approximations are widely used, see, for example, Section 5.1 of \cite{Belotti_Kirches_Leyffer_Linderoth_Luedtke_Mahajan_2013}.  In these cases,~\eqref{sec:surr_problem:tildeP} can be written as MILPs, i.e., they are MILP-representable. In particular, $C$ is the   subset of indices $m\in\{0,\ldots,M\}$ such that the functions $g_m(\mathbf{x})$  are nonlinear. Solving such surrogate problems is particularly interesting as exact methods for  MILPs, on average, perform and scale better than the ones for MINLPs. Some recent papers propose to obtain MILP-representable surrogate problems via data-driven approaches like machine learning techniques, see, e.g., \cite{bertsimas2023global,bertsimas2025global}. In these papers as well, the set $C$  contains the indices of all nonlinear functions in~\eqref{sec:surr_problem:P}.

In this work, we follow a different strategy to build~\eqref{sec:surr_problem:tildeP} from~\eqref{sec:surr_problem:P}. Instead of replacing all nonlinearities in~\eqref{sec:surr_problem:P} with linear approximations to obtain~\eqref{sec:surr_problem:tildeP}, we propose to retain the ``simple'' nonlinearities, i.e., those that are computationally tractable, such as convex,  quadratic or separable univariate terms. Therefore, we focus only on the complicating functions. These are replaced by nonlinear approximations obtained as explained in Section~\ref{sec:shape_constrained_regression}, i.e., piecewise polynomials that might have degree $d> 1$. Then, in our setting,~\eqref{sec:surr_problem:tildeP} is still an MINLP. Our intuition is twofold: on one side, we think the resulting surrogate MINLP corresponds to a better approximation of~\eqref{sec:surr_problem:P}, since the approximations are allowed to be nonlinear; on the other side, we trust that simpler MINLPs, i.e., those where the original complicating nonlinearities are replaced by separable piecewise-polynomials,  could be solved to global optimality by available solvers. In particular, MILP solvers were extended to deal with a quadratic objective and constraints, see, for example, \cite{gurobi,cplex,SCIPOptSuite10}. Moreover, practically effective methods exist for MINLP with sum-of-univariate objective and constraints, see, for example, \cite{d2012algorithmic,d2019strengthening}. The methods described in Section~\ref{sec:shape_constrained_regression} could be used to impose some special properties on the nonlinear approximations of the objective function or constraints in~\eqref{sec:surr_problem:P}, namely, we could impose that the approximation is a pointwise underestimation or overestimation of the original functions. The former would result in problem \eqref{sec:surr_problem:tildeP} being a pointwise relaxation of~\eqref{sec:surr_problem:P} and the latter being a pointwise restriction of~\eqref{sec:surr_problem:P}, i.e., only at the points where the underestimation or overestimation is enforced.

The methods based on solving~\eqref{sec:surr_problem:tildeP} such as it is a restriction or, more in general, an approximation of~\eqref{sec:surr_problem:P},
are heuristic for the following two reasons: 

\begin{enumerate}
\item When we replace the objective function $g_0(\mathbf{x})$ with a function ${\tilde g}_0(\mathbf{x})$, we solve to optimality problem \eqref{sec:surr_problem:tildeP}, i.e., an approximation of problem~\eqref{sec:surr_problem:P}. Clearly, given the optimal solution $\tilde{\mathbf{x}}$ of problem  \eqref{sec:surr_problem:tildeP}, one could always recompute the corresponding original objective function $g_0(\tilde{\mathbf{x}})$. However, the optimal solution of the problem \eqref{sec:surr_problem:tildeP} is, in general, nonoptimal for~\eqref{sec:surr_problem:P}.  The quality of $\tilde{\mathbf{x}}$ with respect to the original objective function cannot be guaranteed and a low-quality solution can be obtained.
\item When we replace (some of) the constraint functions $g_m(\mathbf{x})$ with their surrogate functions ${\tilde g}_m(\mathbf{x})$ (for (some) $m \in  \{1,\dots,M\}$), the optimal solution of the surrogate problem \eqref{sec:surr_problem:tildeP} might be infeasible for~\eqref{sec:surr_problem:P}, i.e., $\exists m: g_m(\tilde{\mathbf{x}}) > 0$. Or, worse, \eqref{sec:surr_problem:tildeP} might have no feasible solution at all for highly constrained  MINLPs. In the case of feasible \eqref{sec:surr_problem:tildeP}, an additional procedure is needed to recover feasibility. 
\end{enumerate}

\subsection{A formulation for the surrogate problem}
\label{sec:surrog_formulation}

 In this section, we introduce a formulation for~\eqref{sec:surr_problem:tildeP} from Definition~\ref{def:surrogate}, where the approximations $\tilde{g}_m(\mathbf{x}),$ $m\in C,$ are constructed as described in  Section~\ref{sec:shape_constrained_regression} and are therefore piecewise univariate polynomials. From now on, we consider $m \in C$ when $g_m(\mathbf{x})$ is nonconvex, $m \not\in C$ when $g_m(\mathbf{x})$ is convex. Consequently, \eqref{sec:surr_problem:tildeP}  can be formulated as an MINLP where the only nonconvex functions can be written as sum of univariate functions.

For simplicity, let us assume that $C = \{0\}$, i.e., our only complicating (or nonconvex) function is $g_0(\mathbf{x})$. For the ease of notation, we drop the index corresponding to the complicating functions, as we assume that we have just one of them. Also, suppose that we consider only the internal knots of each variable $x_j$, i.e., those contained within the observed domain of $x_j$, $[\underline{x}_j,\overline{x}_j]$. From now on, we index the knots as $t_{j1}, t_{j2}, \dots, t_{jk_j}, t_{j(k_j+1)}$, $j=1,\dots, p$,  being $t_{j1}=\underline{x}_j$ and $t_{j(k_j+1)}=\overline{x}_j$.

As we already mentioned, the  function $\tilde{g}(\mathbf{x})$ approximating the complicating objective function $g(\mathbf{x})$ can be written as a sum of univariate piecewise polynomials, formally,

\begin{equation*}
    {\tilde g}(\mathbf{x}) = \sum_{j=1}^p {\tilde g}'_{j}(x_j),
\end{equation*}

\noindent where, for $ j = 1, \dots, p,$

\begin{equation*}
    {\tilde g}'_{j}(x_j) = 
                        \begin{cases}
                        c_{j10} + \sum_{d=1}^{d_j} c_{j1d} (x_{j} - t_{j1})^d & \mbox{for } x_j \in [t_{j1}, t_{j2}] \\
                        c_{j20} + \sum_{d=1}^{d_j} c_{j2d} (x_{j} - t_{j2})^d & \mbox{for } x_j \in [t_{j2}, t_{j3}] \\
                        \dots & \\
                        c_{j({k_j})0} + \sum_{d=1}^{d_j} c_{j({k_j})d} (x_{j} - t_{jk_j})^d & \mbox{for } x_j \in [t_{jk_j}, t_{j({k_j}+1)}].
                        \end{cases}
\end{equation*}

Thus, the surrogate  problem associated with~\eqref{sec:surr_problem:P}, i.e.~\eqref{sec:surr_problem:tildeP}, obtained by approximating   its objective function via a smooth additive regression model, can be formulated as the following MINLP:

\begin{equation}\label{eq:sur_form1}
    \begin{aligned}
& \alpha + \min \sum_{j=1}^p \sigma_j  \\
& \text{s.t.}\\
& \quad g_m(\mathbf{x}) \leq 0, \quad \quad \quad  m \in \{ 1,\dots,M\}, \\
& \quad x_j \in \mathbb{Z},\quad \quad \quad \quad  \quad  j  \in I \subseteq \{1,\dots,p\},\\
& \quad \underline{x}_j \leq x_j \leq \overline{x}_j, \quad \quad   j  \in \{1,\dots,p\},\\
& \quad \mathbf{x} \in \mathbb{R}^p,\\
\end{aligned}
\end{equation}

\noindent where the following additional variables are needed for all $j=~1,\dots,p$ and $q=~1,\dots, k_j$:

\begin{itemize}
    \item $\sigma_j$ represents the value of the univariate surrogate function of variable $x_j$.
    \item $\sigma'_{jq}$ represents the contribution of each interval $[t_{jq},t_{j(q+1)}]$ to the surrogate function of $x_j$. 
    \item $x_{jq}$ is the deviation of $x_j$ from  $t_{jq}$ if  $x_j  \in [t_{jq}, t_{j(q+1)}]$, and zero otherwise.
    \item $y_{jq}$ takes value $1$ if $x_j  \in [t_{jq}, t_{j(q+1)}]$, and zero otherwise. 
\end{itemize}

\noindent The following constraints are also needed:

\begin{align}
&\sigma'_{jq} = c_{jq0} \cdot y_{jq} +\sum_{d=1}^{d_{j}} c_{jqd} \cdot \left(x_{jq} \right)^{d}, & \quad   j=1,\dots,p,  q=1,\dots,k_j, \label{eq:surr_sigma'_def}\\
&0 \leq x_{jq} \leq (t_{j(q+1)} - t_{jq}) \cdot y_{jq}, & \quad   j=1,\dots,p,  q=1,\dots,k_j ,\label{eq:surr_bounds_xjq}\\
&\sum_{q=1}^{k_{j}} y_{jq} = 1, & \quad   j=1,\dots,p, \label{eq:surr_onlyoney} \\
& \sum_{q=1}^{k_{j}} \sigma'_{jq} = \sigma_{j},  & \quad  j=1,\dots,p, \label{eq:sum_sigma'}\\
\end{align}

\begin{align}
&x_{j} = \sum_{q=1}^{k_{j}} \left( y_{jq} \cdot t_{jq} + x_{jq} \right), & \quad   j=1,\dots,p ,\label{eq:surr_xj_def}\\
&y_{jq} \in \{0,1\},  &\quad   j=1,\dots,p,  q=1,\dots,k_j ,\label{eq:surr_binary_y} \\
&\sigma'_{jq}  \in \mathbb{R}, & \quad   j=1,\dots,p,  q=1,\dots,k_j \label{eq:surr_sigma'_bounds}, \\ 
&\sigma_{j}  \in \mathbb{R}, & \quad   j=1,\dots,p\label{eq:surr_sigma_bounds}. 
\end{align}

\noindent  Constraints~\eqref{eq:surr_sigma'_def} and~\eqref{eq:surr_bounds_xjq} model the correct definition of variables $\sigma'_{jq}$ and $x_{jq}$, respectively. Let us note that, for a given $j \in \{1,\dots,p\}$, only one $y_{jq}$ can take value~1 thanks to constraints~\eqref{eq:surr_onlyoney}. Constraints~\eqref{eq:sum_sigma'} ensure that $\sigma_j$ is the sum of  $\sigma'_{jq}$ over all the $k_j$ intervals. Constraints~\eqref{eq:surr_xj_def} links variables $x_j$, $x_{jq}$ and $y_{jq}$. In particular, the summation on the right-hand-side shows only one nonzero term, i.e., the one corresponding to the interval where $x_j$ lies. Let $q$ be the index corresponding to such an interval. Then, $x_j$ is equal to $t_{jq}+x_{jq}$. Constraints~\eqref{eq:surr_binary_y},~\eqref{eq:surr_sigma'_bounds} and~\eqref{eq:surr_sigma_bounds} state the nature of all the variables. Note that variables $\sigma_j$ are summed up in the objective function, together with the intercept constant term $\alpha$. 
 
We note that~\eqref{eq:sur_form1}--\eqref{eq:surr_sigma_bounds} define an MINLP formulation of  problem~\eqref{sec:surr_problem:tildeP} from Definition~\ref{def:surrogate}, with $C=\{0\},$ where the approximation is obtained using the constrained  smooth additive regression model described in Section~\ref{sec:shape_constrained_regression}. We point out that in the case of $C=\{c\}$ where $1\leq c\leq M,$ namely a complicating constraint is being approximated to get~\eqref{sec:surr_problem:tildeP} instead of the objective function, the approach to get the MINLP formulation is analogous.
Finally, let us remark that this mathematical model for~\eqref{sec:surr_problem:tildeP} is inspired by the classical multiple choice formulation, see, for example,  \cite{vielma2009mixed}.

\subsection{The overall algorithm}
\label{sec:surrog_the-overall-algorithm}

 Let us consider that we are given~\eqref{sec:surr_problem:P}, namely a mathematical optimization problem formulated as an MINLP with a non-empty feasible solution set, jointly with an index set $C$ of its complicating  functions  $C \subset \{0,\dots,M\}$. For the sake of readability, we assume that there is only one complicating function, i.e., $C=\{c\}$ is a singleton where $0\leq c\leq M.$ The extension to several complicating functions is straightforward. Note that  the only complicating function, $g_c(\mathbf{x}),$ may appear either in the objective function ($c=0$) or in one of the constraints ($1 \leq c \leq M$).

The MISSOC algorithm consists of four main steps: i) generate a dataset $\mathcal{T}$ from~\eqref{sec:surr_problem:P} by sampling $g_c$; ii) fit an approximating function $\tilde{g}_c$; iii) formulate~\eqref{sec:surr_problem:tildeP} by replacing $g_c$ with $\tilde g_c$ and solve it; and iv) post process the so-obtained solution of~\eqref{sec:surr_problem:tildeP} to obtain a heuristic solution for~\eqref{sec:surr_problem:P}  by recovering feasibility, if needed, and improving the objective value. In what follows, we describe the different steps in detail and present the corresponding pseudocode in Algorithm~\ref{alg:missoc}.

To   formulate a surrogate problem~\eqref{sec:surr_problem:tildeP} 
from~\eqref{sec:surr_problem:P} in the form described in Section~\ref{sec:surrog_formulation}, we first construct a training dataset  by sampling values of $\mathbf{x}$ within the domain defined in~\eqref{sec:surr_problem:P} and computing the corresponding values $g_c(\mathbf{x}).$  These serve as the covariates and response variable, respectively, in the constrained smooth additive regression model used to approximate the complicating function. Then, let this training data set $\mathcal{T}$ consist of $n$ samples
of the form $\mathcal{T} =~\{(\mathbf{x}_i, g_c(\mathbf{x}_i)),\, i=~1,\ldots,n\}$. Due to our sampling procedure, we need to assume that all  the decision variables $x_j$ appearing in the function $g_c(\mathbf{x})$ have finite bounds, i.e., $\underline{x}_j \neq -\infty$ and $\overline{x}_j \neq +\infty$. Ideally, 
 the sampled $\mathbf{x}_i,$ $i=1,\ldots,n,$ in $\mathcal{T}$   should be well spread and space-filling over the bounded domain of interest in order to provide a representative description of the function $g_c$. In MISSOC, we use uniform random sampling to obtain $\mathcal{T}$, which in general  provides a good compromise between computational cost and the quality of the training data. However, other sampling strategies are possible but are out of the scope of this paper. The reader is referred to \cite{garud2017design} for a recent review on this topic.

Once the training data has been generated, the function $\tilde{g}_c$, approximating $g_c$, is obtained by means of a smooth additive regression model with constraints, as described in Section~\ref{sec:shape_constrained_regression}. The function is estimated using training data $\mathcal{T}$, and a selected degree $d_j$ of the $B-$spline basis, and number of intervals $k_j$ for each of the $p$  decision variables acting as covariates for the regression model.

The third step consists of constructing the surrogate problem~\eqref{sec:surr_problem:tildeP} using the previously estimated function $\tilde{g}_c$ and its MINLP formulation, as explained in Section~\ref{sec:surrog_formulation}. 
 Observe that in the case of $c=0$,~\eqref{sec:surr_problem:tildeP} is always feasible, while this is not necessarily true if $1 \leq c \leq M$. To address this feasibility issue arising when approximating a constraint, we consider the following relaxation of~\eqref{sec:surr_problem:tildeP}, in which a slack variable $u$ is introduced into the approximated constraint:

\begin{equation}
\begin{aligned}
& \min \quad g_0(\mathbf{x}) \\
& \text{ s.t.} \\
& \quad \tilde g_c(\mathbf{x}) - u \leq 0, \\
& \quad g_m(\mathbf{x}) \leq 0, \quad m \in \{1,\dots,M\},\ \ m \neq c, \\
& \quad x_j \in \mathbb{Z}, \quad j \in I \subseteq \{1,\dots,p\}, \\
& \quad \underline{x}_j \leq x_j \leq \overline{x}_j, \quad j \in \{1,\dots,p\}, \\
& \quad \mathbf{x} \in \mathbb{R}^p, \\
& \quad 0 \leq u \leq \overline{u}.
\end{aligned}
\tag{$\tilde{P}_R$}
\label{alg:relaxation}
\end{equation}

\noindent Here, the upper bound $\overline{u}$ controls the maximum violation allowed for the approximated constraint in the relaxed~\eqref{sec:surr_problem:tildeP}, namely~\eqref{alg:relaxation}. To define $\overline{u}$, we use the maximum absolute approximation error $\hat \varepsilon$ observed on the training dataset $\mathcal{T}$ given by

\begin{equation*}
   \hat \varepsilon= \max_{i=1,\dots,n} \left| g_c(\mathbf{x}_i)-\tilde g_c(\mathbf{x}_i)\right|.
\end{equation*}

\noindent Then, we choose $\overline{u}$ as the smallest power of 10 that is greater than or equal to $\hat \varepsilon$, providing a simple and conservative tolerance that reflects the observed pointwise approximation error.

The fourth step in MISSOC consists of solving the surrogate problem~\eqref{sec:surr_problem:tildeP} (resp.~\eqref{alg:relaxation}) Recall that $\tilde{g}_c$ is an additive regression model, which consists of the sum of univariate polynomials. In other words,~\eqref{sec:surr_problem:tildeP}  (resp.~\eqref{alg:relaxation}) is formulated as an MINLP in which nonconvexities appear as sum of univariate functions. Thus, we propose to exploit this separable structure and solve it with the SC-MINLP algorithm \citep{d2012algorithmic,d2019strengthening}, which is specifically designed for MINLPs with separable (sum of univariate) nonconvex functions and is therefore particularly suitable for solving the surrogate problem generated by MISSOC. Nevertheless, any global optimization solver capable of handling mixed-integer nonlinear polynomial functions could be used to solve~\eqref{sec:surr_problem:tildeP}  (resp.~\eqref{alg:relaxation}) at the expense of not specifically exploiting its structure.

 Provided that $\tilde{\mathbf{x}}$ is the solution obtained in the previous step, two situations appear. 
If the approximated function is the objective function ($c=0$), $\tilde{\mathbf{x}}$  is obtained by solving~\eqref{sec:surr_problem:tildeP}. Therefore, $\tilde{\mathbf{x}}$ is feasible in the original problem~\eqref{sec:surr_problem:P} and it is a heuristic solution for it, with objective value $g_0(\tilde{\mathbf{x}})$.  Nevertheless, $\tilde{\mathbf{x}}$ might not be optimal for~\eqref{sec:surr_problem:P} and its quality cannot be guaranteed. If, instead, the approximated function appears in a constraint ($1 \leq c \leq M$), then $\tilde{\mathbf{x}}$ is obtained by  solving ~\eqref{alg:relaxation} and it may be infeasible for~\eqref{sec:surr_problem:P}. Therefore, $\tilde{\mathbf{x}}$ is not, in general, a heuristic solution for the original problem. To handle both situations, the last step of MISSOC is a post processing phase that attempts to recover feasibility with respect to~\eqref{sec:surr_problem:P}, when necessary, and to improve the quality of the obtained solution. In the following, we denote by $\mathbf{x}^*$ to the final heuristic feasible solution for the original problem~\eqref{sec:surr_problem:P} returned by MISSOC.

\begin{algorithm} 
\caption{\textbf{MISSOC}\label{alg:missoc}}
\textbf{Input:} Problem \eqref{sec:surr_problem:P}. \\
\hspace*{3.2em}  Index of complicating function $c$. \\
\hspace*{3.2em} Degrees $(d_1,\ldots,d_p)$. \\ 
\hspace*{3.2em} Number of intervals $(k_1,\ldots,k_p)$.\\
\hspace*{3.2em} Number of samples $n$\\
\textbf{Output:} Global optimal solution $\tilde{\mathbf{x}}$ for~\eqref{sec:surr_problem:tildeP}. \\
\hspace*{4.1em} Heuristic solution $\mathbf{x}^*$ to~\eqref{sec:surr_problem:P}.
\begin{algorithmic}[1]
    \State Using random sampling, generate a training dataset $\mathcal{T} =~\{(\mathbf{x}_i, g_c(\mathbf{x}_i)),\, i=~1,\ldots,n\}$.
    \State $\tilde{g}_c \gets $ estimate $g_c$ as in Section~\ref{sec:shape_constrained_regression} using data in $\mathcal{T}$,  degrees $(d_1,\ldots,d_p)$, and number of intervals $(k_1,\ldots,k_p).$
    \State Formulate~\eqref{sec:surr_problem:tildeP} (resp.~\eqref{alg:relaxation}) as in Section~\ref{sec:surrog_formulation} using $\tilde{g}_c$.
    \State $\tilde{\mathbf{x}} \gets$ solution of~\eqref{sec:surr_problem:tildeP} (resp.~\eqref{alg:relaxation}) with a global optimization solver.
    \State $\mathbf{x}^*\gets$  apply the post processing phase (Algorithm~\ref{alg:postprocessing}) to $\tilde{\mathbf{x}}$.\\
   \Return $\mathbf{x}^*$.
\end{algorithmic}
\end{algorithm}

The post processing phase performs a local search around the solution $\tilde{\mathbf{x}}$ through three consecutive subphases. First,  the values of the integer variables are fixed and only the continuous variables are allowed to change. This may improve the objective value and potentially recover feasibility. In a second subphase, if feasibility  has not been recovered, the local search is progressively enlarged by allowing limited changes in the integer variables through a local branching constraint \citep{fischetti2003local}. Finally, after a first feasible solution is found, an optional third subphase further enlarges these local branching neighborhoods to improve the objective value while maintaining feasibility. Note that this third subphase is optional since the feasible solution obtained in the second subphase may already be sufficient in some applications where quickly obtaining a feasible solution is prioritized over further objective improvements. Each
of these three subphases is detailed below and the complete post processing step in
MISSOC (step 5) is summarized in Algorithm~\ref{alg:postprocessing}.

The local search in the post processing is controlled by a local branching constraint that limits the distance between the integer components of  a new solution $\mathbf{x}$ and those of  the incumbent $\tilde{\mathbf{x}}$. In general, this can be written as
\begin{equation}
    D_I(\mathbf{x},\tilde{\mathbf{x}}) \leq \kappa,
\label{alg:MISSOC_localbranching_general}
\end{equation}
\noindent where $D_I$ denotes a distance restricted to the integer variables and $\kappa$ controls the size of the search neighborhood.  For simplicity, in the description of the post processing step we assume that the integer variables are binary. In this case, $D_I$ is the Hamming distance and $\kappa$ has the direct interpretation of the maximum number of binary variables in $\mathbf{x}$ whose value may differ from that of $\tilde{\mathbf{x}}$. Hence,~\eqref{alg:MISSOC_localbranching_general} can be written as the linear constraint
\begin{equation}
\sum_{\substack{j \in I\\ \tilde{x}_j=0}} x_j  
+ 
\sum_{\substack{j \in I\\ \tilde{x}_j=1}} (1-x_j) 
\leq \kappa.
\label{alg:MISSOC_localbranching}
\end{equation}

\noindent For general integer variables, other formulations based on distances can be used, at the expense of introducing additional variables and constraints.

The problem considered in the three subphases of the MISSOC post processing step is described through the formulation of an MINLP, which includes the local branching constraint~\eqref{alg:MISSOC_localbranching} to define the binary search neighborhood, and stated as:

\begin{equation}
\begin{aligned}
& \min \quad \omega_0 g_0(\mathbf{x}) + \omega_{u}u/\overline{u} \\
& \text{ s.t.} \\
& \quad g_c(\mathbf{x}) - u \leq 0,\\
& \quad g_m(\mathbf{x}) \leq 0, \quad m \in \{ 1,\dots,M\},\ \ m \neq c,  \\
&\quad \sum_{\substack{j \in I\\ \tilde{x}_j=0}} x_j + \sum_{\substack{j \in I\\ \tilde{x}_j=1}} (1-x_j) \leq \kappa, \\
& \quad x_j \in \{0,1\}, \quad j \in I \subseteq \{1,\dots,p\}, \\
& \quad \underline{x}_j \leq x_j \leq \overline{x}_j, \quad j \in \{1,\dots,p\}, \\
& \quad \mathbf{x} \in \mathbb{R}^p, \\
& \quad 0 \leq u \leq \overline{u}.
\end{aligned}
\label{alg:ls_problem}
\end{equation}

\noindent Problem~\eqref{alg:ls_problem} is built from the original problem~\eqref{sec:surr_problem:P}. It is inspired by custom feasibility restoration strategies, where constraint violation is modeled and reduced through auxiliary slack variables. The original complicating constraint $g_c$ is relaxed through the slack variable $ u\geq 0$, upper bounded by $\overline{u}$. The non-complicating constraints of~\eqref{sec:surr_problem:P}, i.e., $g_m$ such that  $m \in \{ 1,\dots,M\}$ and $m \neq c$, remain unchanged without relaxation. The objective function in~\eqref{alg:ls_problem} combines the objective function in~\eqref{sec:surr_problem:P}, $g_0$, and the normalized feasibility violation term $u/\overline{u}$, with weights $\omega_0$ and $\omega_u$, respectively. The weights satisfy $\omega_0+\omega_u=1$, where  $\omega_0,\omega_u\in[0,1]$, and control the trade-off between objective improvement and feasibility recovery. For the two terms in the objective of problem~\eqref{alg:ls_problem} to be comparable,  we assume that $g_0$ is scaled to take values in $[0,1]$ or that it has been normalized accordingly. Problem~\eqref{alg:ls_problem} also includes the local branching constraint~\eqref{alg:MISSOC_localbranching}, which  limits to $\kappa$ the number of binary variables allowed to change with respect to $\tilde{\mathbf{x}}$. As discussed above, the formulation is stated for binary integer variables, so that~\eqref{alg:MISSOC_localbranching} can be applied.

 We point out that problem~\eqref{alg:ls_problem} is stated for the case in which the approximated function is a constraint, i.e., $1 \leq c \leq M$. When $c=0$, the same formulation holds but omitting the slack variable $u$, the relaxed constraint involving it and the corresponding term in the objective. One could simply obtain it by setting $\overline{u}=0$. In what follows, we describe the post processing step for the case  $1 \leq c \leq M$. The case $c=0$ follows directly by omitting the terms involving $u$, i.e., by setting $\overline{u}=0$.

Problem~\eqref{alg:ls_problem} is used throughout the three subphases of the post processing step to restore feasibility or/and improve the quality of $\tilde{\mathbf{x}}$. This problem is solved iteratively, starting from $\kappa=0$ and progressively increasing it to enlarge the local branching neighborhood around $\tilde{\mathbf{x}}$. The weights $\omega_0$ and $\omega_u$ and the slack upper bound $\overline{u}$ in problem~\eqref{alg:ls_problem} are parameters that are set depending on the subphase considered. We denote by $\mathbf{x}^{(\kappa)}$ and $u^{(\kappa)}$ the solution of problem~\eqref{alg:ls_problem} at iteration $\kappa$, when it exists.

In the first subphase of the post processing step, problem~\eqref{alg:ls_problem} is solved by setting $\kappa=0$ and using the same slack upper bound $\overline{u}$ as in the relaxed surrogate problem~\eqref{alg:relaxation}, or $0$ if $c=0$. The weights $\omega_0$ and $\omega_u$ are chosen by the user, who may prioritize feasibility recovery, objective improvement or a combination of both. Since $\kappa=0$, the local branching constraint fixes the binary variables to their values in $\tilde{\mathbf{x}}$. Therefore, this first subphase reduces to a relaxed restricted nonlinear program (NLP) of the original problem~\eqref{sec:surr_problem:P}, where only the continuous variables are allowed to change and the violation of the complicating constraint is bounded by $\overline{u}$. The resulting problem is nonconvex and can be solved to local optimality using $\tilde{\mathbf{x}}$ as a warm start solution. As a result, this first subphase may improve the objective value and may also recover feasibility with respect to~\eqref{sec:surr_problem:P}.
If problem~\eqref{alg:ls_problem} is feasible and $u^{(0)}$ is smaller than or equal to a fixed feasibility tolerance, which we set to $10^{-6}$ following standard practices, then $\mathbf{x}^{(0)}$ is considered feasible for~\eqref{sec:surr_problem:P} and is stored as the current best feasible solution, denoted by $\mathbf{x}^{\mathrm{best}}$. If no further improvement is desired, $\mathbf{x}^{\mathrm{best}}$ is returned as the heuristic solution $\mathbf{x}^*$ of~\eqref{sec:surr_problem:P}. Otherwise, the post processing may proceed to the optional third subphase for objective improvement.

Note that, in the first subphase, feasibility recovery with respect to~\eqref{sec:surr_problem:P} is not guaranteed: problem~\eqref{alg:ls_problem} may either have no solution or have one with $u^{(0)}>10^{-6}$. When binary variables are present, feasibility may require changes in their values.  In such cases, the post processing proceeds to the second subphase.

In the second subphase, feasibility recovery is addressed by allowing local changes in binary variables with respect to the surrogate solution $\tilde{\mathbf{x}}$. To this end, the weights in problem~\eqref{alg:ls_problem} are set to $\omega_0=0$ and $\omega_u=1$, so that the objective minimizes only the feasibility violation $u$. The second subphase starts by solving problem~\eqref{alg:ls_problem} with $\kappa=1$ and with $\overline{u}$ set as in~\eqref{alg:relaxation}. If $\mathbf{x}^{(0)}$ exists, it is used as a warm start; otherwise, $\tilde{\mathbf{x}}$ is used. The process then continues iteratively by increasing $\kappa$ by one to enlarge the local branching neighborhood around $\tilde{\mathbf{x}}$. At each iteration $\kappa$,  the problem is warm-started with the solution from the previous iteration, $\mathbf{x}^{(\kappa-1)}$ and $\overline{u}$ is set to $u^{(\kappa-1)}$. In the case where a solution to problem~\eqref{alg:ls_problem} has not been found yet in previous iterations,  $\overline{u}$  remains unchanged and the warm start is $\tilde{\mathbf{x}}$. The process stops when $u^{(\kappa)} \leq 10^{-6}$, i.e, a feasible solution of~\eqref{sec:surr_problem:P} is found. Let $\kappa_F$ be the value of $\kappa$ at which feasibility is recovered. Then, $\mathbf{x}^{(\kappa_F)}$ is feasible for~\eqref{sec:surr_problem:P} and is stored as the  best solution $\mathbf{x}^{\mathrm{best}}$. Note again that feasibility recovery is not guaranteed in general. However, the underlying expectation is that the surrogate solution $\tilde{\mathbf{x}}$ lies close to a feasible solution of the original problem~\eqref{sec:surr_problem:P}, so that increasing $\kappa$ to explore larger local branching neighborhoods may help recover feasibility.

\begin{algorithm}
\footnotesize
\caption{\textbf{Post processing phase in MISSOC}\label{alg:postprocessing}}

\textbf{Input:} Original problem~\eqref{sec:surr_problem:P}. \\
\hspace*{3.2em}  Index of complicating function $c$.\\
\hspace*{3.2em} Surrogate solution $\tilde{\mathbf{x}}$.\\
\hspace*{3.2em} Initial maximum slack value $ \overline{u}$, if a constraint is approximated.\\
\hspace*{3.2em} Weights $\omega_0,\omega_u$ for the first subphase.\\
\hspace*{3.2em} Feasibility tolerance $\Sigma$.\\
\hspace*{3.2em} Boolean parameter \texttt{Improve}.\\
\hspace*{3.2em} Relative change tolerance $\tau$.\\
\hspace*{3.2em} Patience parameter $\eta$.\\
\textbf{Output:} Best solution  found for~\eqref{sec:surr_problem:P}, $\mathbf{x}^{\mathrm{best}} $.

\begin{algorithmic}[1]

\State $\mathbf{x}^{\mathrm{best}} \gets \emptyset$.
\State  $\kappa \gets 0$.
\State Solve problem~\eqref{alg:ls_problem}  using $\tilde{\mathbf{x}}$ as a warm start.
\If{a solution exists}
    \State $\mathbf{x}^{(0)},  \ u^{(0)} \gets $ solution of problem~\eqref{alg:ls_problem}.
\EndIf 
\If{ $\mathbf{x}^{(0)}$ exists and  $u^{(0)}\leq \Sigma$}
\State  $\mathbf{x}^{\mathrm{best}} \gets\mathbf{x}^{(0)}$
\State  $\overline{u} \gets u^{(0)}$.

\Else

    \State  $\kappa \gets 1$.
    \State  $\omega_0 \gets 0$ and $\omega_u \gets 1$.
    \State  \texttt{feasible} $\gets \mathrm{false}$.

    \While{not \texttt{feasible}}

    \State Solve problem~\eqref{alg:ls_problem} using  $\mathbf{x}^{(\kappa-1)}$ as a warm start, if available, and $\tilde{\mathbf{x}}$ otherwise.
    \If{a solution exists}
        \State $\mathbf{x}^{(\kappa)},  \ u^{(\kappa)} \gets $ solution of problem~\eqref{alg:ls_problem}.
        \State $\overline{u} \gets u^{(\kappa)}$.
    \EndIf
        \If{solution $\mathbf{x}^{(\kappa)}$  exists and  $u^{(\kappa)}\leq \Sigma$}
            \State  $\mathbf{x}^{\mathrm{best}} \gets $ $\mathbf{x}^{(\kappa)}$.
            \State \texttt{feasible} $\gets \mathrm{true}$.
        \Else
            \State Set $\kappa \gets \kappa + 1$.
        \EndIf
    \EndWhile
\EndIf

\If{ not \texttt{Improve}}
\State \Return $\mathbf{x}^{\mathrm{best}}$.
\EndIf
    \State  $\omega_0 \gets 1$ and $\omega_u \gets 0$.
    \State  $\mathrm{count} \gets 0$.\\
    \State  $\overline{u} \gets \overline{u} + \Sigma^2$.

 \While{$\mathrm{count} < \eta$}
        \State $\kappa \gets \kappa + 1$.
        \State   $\mathbf{x}^{(\kappa)} \gets$ solution of~\eqref{alg:ls_problem} using $\mathbf{x}^{\mathrm{best}}$ as a starting point.

        \State Compute $
       \Delta^{(\kappa)}$ as in~\eqref{MISSOC:Delta_k}.
             \If{$\Delta^{(\kappa)} < \tau$}
            \State  $\mathrm{count} \gets \mathrm{count}+1$.
        \Else
            \State  $\mathrm{count} \gets 0$.
        \EndIf
          \If{$g_0(\mathbf{x}^{(\kappa)}) < g_0(\mathbf{x}^{\mathrm{best}})$}
            \State $\mathbf{x}^{\mathrm{best}} \gets \mathbf{x}^{(\kappa)}$.
        \EndIf
    \EndWhile
    \State \Return $\mathbf{x}^{\mathrm{best}}$.
\end{algorithmic}

\end{algorithm}

If no further improvement is desired, $\mathbf{x}^{\mathrm{best}}$ is returned as the heuristic solution $\mathbf{x}^*$ of~\eqref{sec:surr_problem:P}. In contrast, if further improvement is desired, the post processing algorithm proceeds to the optional third subphase, which is devoted to improving $g_0$ while maintaining feasibility. To this end, problem~\eqref{alg:ls_problem} is solved with $\omega_0=1$, $\omega_u=0$ and $\overline{u}=u^{(\kappa_F)}+10^{-7}$, where the small tolerance is included to account for numerical errors. The third subphase  solves problem~\eqref{alg:ls_problem} iteratively starting with $\kappa=\kappa_F+1$ and progressively increasing $\kappa$ by one
to enlarge the local branching neighborhood around $\tilde{\mathbf{x}}$. At each iteration $\kappa$, problem~\eqref{alg:ls_problem} is warm-started with the best solution found so far, namely, $\mathbf{x}^{\mathrm{best}}$. 
For the obtained solution $\mathbf{x}^{(\kappa)}$,  we first compute the relative change of its objective value with respect to the best objective value found so far as
\begin{equation}
     \Delta^{(\kappa)}
     =
     \frac{
     g_0(\mathbf{x}^{\mathrm{best}}) - g_0(\mathbf{x}^{(\kappa)})}
     {\left|g_0(\mathbf{x}^{(\kappa)})\right|}.
     \label{MISSOC:Delta_k}
\end{equation}
\noindent Then, if $g_0(\mathbf{x}^{(\kappa)}) < g_0(\mathbf{x}^{\mathrm{best}})$, we update $\mathbf{x}^{\mathrm{best}}=\mathbf{x}^{(\kappa)}$. The process stops when $\Delta^{(\kappa)} < \tau$ for $\eta$ consecutive values of $\kappa$, where $\tau$ is the relative change tolerance and $\eta$ is the patience parameter. At the end of this third subphase, the best solution found $\mathbf{x}^{\mathrm{best}}$ is returned. Algorithm~\ref{alg:missoc} returns it as the heuristic solution $\mathbf{x}^*$ of~\eqref{sec:surr_problem:P}.

The extent to which the three subphases are needed depends on the role of the approximated function. If only the objective function is approximated in MISSOC ($c=0$), the surrogate solution $\tilde{\mathbf{x}}$ is already feasible for~\eqref{sec:surr_problem:P}, regardless of whether integer variables are present or not. In this case, the post processing phase is mainly intended to improve the quality of the surrogate solution and the continuous refinement performed in the first subphase may already be sufficient. However, if there are integers in~\eqref{sec:surr_problem:P}, the third subphase can still be applied if desired to further improve the objective value by allowing controlled changes in their values.

In contrast, if a constraint is approximated ($1 \leq c \leq M$), the surrogate solution $\tilde{\mathbf{x}}$ may be infeasible for~\eqref{sec:surr_problem:P} and the post processing phase is also intended to recover feasibility. The first subphase may recover feasibility when changes in the continuous variables are sufficient, although this is not guaranteed because the corresponding continuous problem is solved only locally. If feasibility is not recovered and integer variables are present, the algorithm proceeds to the second subphase, where local branching allows controlled changes in the integer variables. If feasibility has been recovered, the third subphase may optionally be used to improve the objective value while maintaining feasibility.


\section{Demonstrative example with a real-world problem}
\label{sec:demonstrative}

In this section, we illustrate how MISSOC works through the well-known Water Distribution Network problem \citep{bragalli2012optimal}. It is a real-world problem that consists of designing a water distribution network. The network is composed of nodes connected by pipes. 
The nodes represent points where water is supplied, demanded, or transferred between pipes. Each pipe must be assigned a diameter chosen from a discrete set. The problem contains binary variables that indicate which diameter is selected for each pipe. Since each diameter option has an associated unit cost, the objective is to minimize the total cost of the selected pipe diameters. However, the resulting design must satisfy hydraulic requirements at the nodes to be feasible.

Several instances of the Water Distribution Network problem have been proposed in the literature, differing in the number of pipes and nodes, the topology, and the hydraulic parameters. We focus on the Hanoi instance, with 32 nodes and 34 pipes.

We consider this problem of special interest to showcase MISSOC for three main reasons. First, it is a real-world problem containing integer (binary) variables. Second, the complicating function appears in a constraint and, as shown later, the surrogate solution $\tilde{\mathbf{x}}$  is infeasible for the original problem. Thus, the post processing Algorithm~\ref{alg:postprocessing} has to be run to recover feasibility. Therefore, this example allows us to illustrate the complete MISSOC workflow, including the three post processing subphases. Third, as explained later in Section~\ref{sec:demonstrative:sampling_approx}, the complicating constraint to be approximated satisfies known shape properties, which are incorporated into the approximation. This showcases the ability of MISSOC to build surrogates that are both data-driven and knowledge-driven, which differentiates it from other state-of-the-art surrogate approaches.

The complicating constraint in the Water Distribution Network problem is the hydraulic head loss imposed on each pipe $e$ in the set of pipes $E$. This constraint involves a nonlinear term of the form
\begin{equation}
     \operatorname{sgn}(Q(e)) |Q(e)|^{1.852} A(e)^{-2.435},
     \label{demonstrative:WDN_constraint}
\end{equation}
\noindent where $Q(e)$ denotes the flow through pipe $e$ and $A(e)$ is the  area of the pipe. We refer the reader to \citet{bragalli2012optimal} for a complete description of the problem. The term in~\eqref{demonstrative:WDN_constraint} involves non-integer powers, which can create numerical difficulties for general-purpose solvers.

Since pipe diameters are selected from a finite set of available options, we use an extended formulation in which the head loss constraint is imposed separately for each pipe and for 
each diameter option. In this formulation, the pipe area is thus fixed for each diameter option, so that the term $A(e)^{-2.435}$  in~\eqref{demonstrative:WDN_constraint} becomes a constant. 

In the Hanoi instance, there are six possible diameters. The extended formulation leads to a problem with 712 continuous variables, 204 binary variables and 1697 constraints. Since the network has 34 pipes and each pipe can take one of the 6 possible diameters, $34 \cdot 6 = 204$ of these constraints are the complicated ones that contain the nonlinear 
term~\eqref{demonstrative:WDN_constraint}.

{
\renewcommand{\arraystretch}{1.2}
\begin{table}[h!]
 \centering
\begin{tabular}{@{}crrr@{}}
\toprule

        \textbf{Metric} & \multicolumn{1}{c}{\textbf{Gurobi}}        
         & \multicolumn{1}{c}{\textbf{BARON}} & \multicolumn{1}{c}{\textbf{SCIP}}\\ \toprule
         $g_0(\mathbf{x}^{(P)})$   & 6109620.900&   n/a &  6109620.900 \\
         gap(\%) &0& n/a  & 0     \\
         time(s)  & 386.370  & 600  & 576.670    \\ 
         \bottomrule
\end{tabular}
\caption{Results for the original problem~\eqref{sec:surr_problem:P} of the Hanoi instance with the extended formulation.}
\label{tab:WDN_hanoi_original}
\end{table}
}

Table~\ref{tab:WDN_hanoi_original} reports the results obtained by solving the original extended formulation of the Hanoi instance with Gurobi v.\,12.0.1 \citep{gurobi}, BARON v.\,24.5.8 \citep{sahinidis1996baron}, and SCIP v.\,10.0.0 \citep{SCIPOptSuite10}. We included SCIP because it was extended to deal with the nonlinearities arising from the Water Distribution Network problem, see \cite{vigerske2013decomposition}.  All solvers were run with a time limit of 600 seconds and with the rest of the parameters set to default. For each solver, we report the objective value achieved,  $g_0(\mathbf{x}^{(P)})$, being $\mathbf{x}^{(P)}$ the corresponding solution found. We also report the optimality gap reported by the solver and the computational time in seconds. We also attempted to solve the instance with COUENNE 0.5.8 \citep{belotti2009branching} and BONMIN \citep{bonami2008algorithmic}, but these solvers returned errors when handling non-integer power terms. Gurobi and SCIP are able to solve the extended formulation to optimality. Gurobi closes the optimality gap in 386.37 seconds, although the optimal solution is found after 351 seconds. SCIP closes the gap in 576.67 seconds, finding the optimal solution after 475 seconds. In contrast, BARON reaches the time limit without finding a feasible solution. For completeness, we also tested the non-extended formulation of the Hanoi instance, for which the solvers showed worse performance.

In what follows, we apply MISSOC step by step to the Hanoi instance, showing how the complicating constraints are approximated, how the surrogate problem is constructed and solved and, finally, how the post processing phase recovers feasibility and improves the final solution.

\subsection{Sampling and constraint approximation}
\label{sec:demonstrative:sampling_approx}

Recall that the complicating constraints to approximate in the original problem are the hydraulic head loss relations, which contain the nonlinear term in~\eqref{demonstrative:WDN_constraint}. Since we use the extended formulation, this head loss constraint is imposed separately for each pipe and each diameter option. Let $\mathcal{D}=\{1,\ldots,6\}$ 
denote the index set of available diameters in the Hanoi instance and let $\delta_r, \ r \in \mathcal{D},$ be the $r$-th diameter. Let $a_r$ be the pipe area associated with diameter $\delta_r$. Then, for each fixed diameter $\delta_r, \ r \in \mathcal{D},$  the area  $A(e)$ in~\eqref{demonstrative:WDN_constraint} becomes the constant $a_r$, so the nonlinear term to be approximated is
\begin{equation}
    \psi(Q(e)) = \operatorname{sgn}(Q(e))|Q(e)|^{1.852}.
   \label{demonstrative:term_approximated}
\end{equation}

\noindent Although the analytical form of~\eqref{demonstrative:term_approximated} is the same for all diameters, the range of $Q(e)$ depends on the selected diameter. In particular, for each $r\in\mathcal{D}$,
\begin{equation}
    -v_{\max}a_r \leq Q(e) \leq v_{\max}a_r, \quad e\in E,
\label{demonstrative:flow_bounds}
\end{equation}
where $v_{\max}$ is also a constant.  We refer again to \citet{bragalli2012optimal} for a detailed description of the full formulation.

Therefore, each diameter $\delta_r, \ r \in \mathcal{D},$ leads to a different domain~\eqref{demonstrative:flow_bounds} for the same nonlinear term~\eqref{demonstrative:term_approximated}. For this reason, we approximate~\eqref{demonstrative:term_approximated} separately over each of the six diameter-dependent ranges of $Q(e)$ in~\eqref{demonstrative:flow_bounds}. Note that these approximations are not pipe-specific. Thus, in the rest of this subsection, we drop the pipe index and write $Q$ for a generic flow value. As shown in Section~\ref{sec:demonstrative:surrogate}, the resulting approximations are later reused across all pipes when building the surrogate problem.

The six approximations of the nonlinear term in~\eqref{demonstrative:term_approximated}, one for each diameter option, are obtained using the shape-constrained smooth additive regression models described in Section~\ref{sec:shape_constrained_regression}. To this end, for each $r\in\mathcal{D}$, we first generate a training dataset of the form $\mathcal{T}^r   = \left\{ \left(Q_i^r, \psi(Q_i^r) \right), \ i=1,\ldots,n     \right\}$, where $Q_i^r$ is sampled uniformly over the
range defined in~\eqref{demonstrative:flow_bounds}.

For each of the six approximations, we use degree $d=3$ and $k=10$ intervals, which leads to a total of $d+k+1=14$ regression coefficients and intercept in the model. We generate 15 observations for each of them, resulting in $n=210$ samples for each training set $\mathcal{T}^{r}$. The rationale behind this sampling rule and the choice of $d$ and $k$ is further discussed later in the experimental setup in Section~\ref{sec:exp:setup}.

\begin{figure}[h!]
    \centering
    \includegraphics[width=0.6\textwidth]{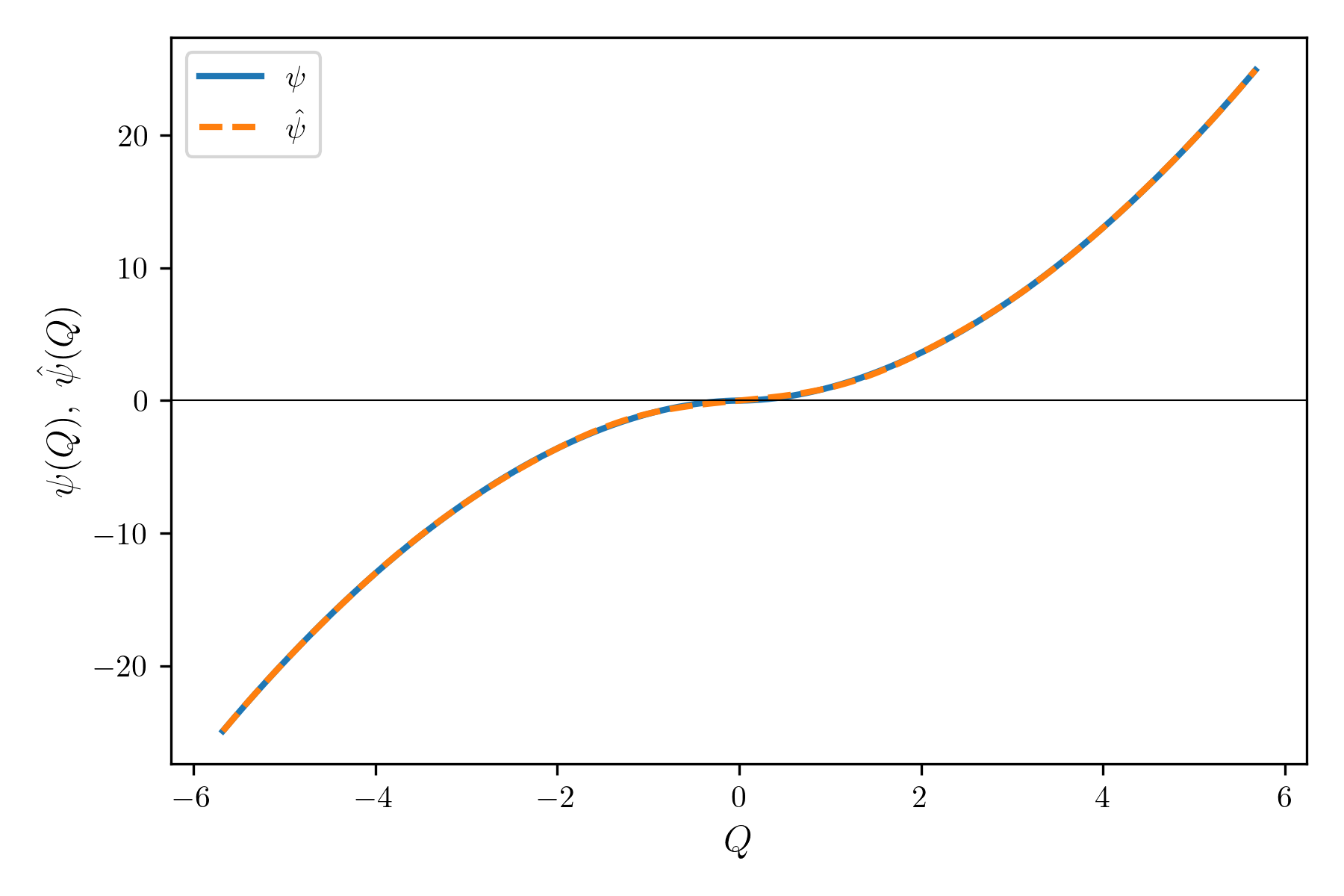}
    \caption{$\psi(Q)$ for the $\delta_6$ diameter  in the Hanoi instance and its approximation $\hat\psi(Q)$.}
    \label{fig:hanoi_approximation}
\end{figure}

Once the training datasets $\mathcal{T}^{r}$ have been generated, we estimate the function $\psi$ in~\eqref{demonstrative:term_approximated} for each diameter $\delta_r, \ r \in\mathcal{D}$. Note that $\psi$ is monotonically nondecreasing with respect to $Q$. Therefore, following the methodology explained in Section~\ref{sec:shape_constrained_regression}, this prior knowledge is incorporated into the estimation of the approximating function $\hat{\psi}$ to ensure that it aligns with the known behavior. In addition, the extended formulation of the Water Distribution Network problem requires that $\psi(0)=0$, so we also impose this pointwise interpolation constraint. 
The total computational time to estimate  $\hat{\psi}$ for the six diameters is 0.268 seconds. 
As an example, Figure~\ref{fig:hanoi_approximation} shows the true function $\psi$ and its approximation $\hat{\psi}$ for the $\delta_6$ diameter in the Hanoi instance. The approximation closely matches the true function over the $Q$ domain and the same behavior is observed for the remaining diameters. The maximum absolute approximation error observed on the training set across all diameters is $\hat{\varepsilon}=0.100$, which occurs in the neighborhood of $Q=0$.

\subsection{Building and solving the surrogate problem}
\label{sec:demonstrative:surrogate}

Once the six approximations of $\psi$ have been obtained, one for each possible diameter, the surrogate problem~\eqref{sec:surr_problem:tildeP} is built by replacing each nonlinear term $\psi$ with its approximation $\hat \psi$, as described in Section~\ref{sec:surrog_formulation}. Recall that, in the extended formulation, the complicating hydraulic head loss constraint is imposed for each pipe and for each diameter. Therefore, for each pair $(e,r)\in E\times\mathcal{D}$, the original term $\psi(Q(e))$ is replaced by the fitted approximation corresponding to the $r$-th diameter. Thus, the six fitted approximations are reused across all 34 pipes.

Since the complicating functions in~\eqref{sec:surr_problem:P} appear in the constraints, the surrogate problem~\eqref{sec:surr_problem:tildeP} is relaxed as in~\eqref{alg:relaxation}. As reported in Section~\ref{sec:demonstrative:sampling_approx}, the maximum absolute approximation error observed on the training sets across all diameters is $\hat{\varepsilon}=0.100$. Therefore, the slack upper bound $\overline{u}$ in~\eqref{alg:relaxation} is set to 1, following the policy stated in Section~\ref{sec:surrog_the-overall-algorithm}.

The resulting relaxed surrogate problem~\eqref{alg:relaxation} has 3161 continuous variables, 2244 binary variables, and 4757 constraints. As explained in Section~\ref{sec:surrog_formulation}, this increase in size with respect to the original problem is due to the multiple choice formulation used to embed the fitted approximations, which requires additional binary variables, continuous variables, and constraints.

We solve~\eqref{alg:relaxation} using the same solvers considered for the original problem with a time limit of 600 seconds and default parameter settings. The results are reported in Table~\ref{tab:WDN_hanoi_surrogate}, which contains the objective value at the surrogate solution, i.e., $g_0(\tilde{\mathbf{x}})$, the optimality gap reported by each solver and the computational time in seconds. 
Gurobi and BARON solve~\eqref{alg:relaxation} to optimality, while SCIP cannot find a feasible solution. Since, in this case, the constraints of~\eqref{sec:surr_problem:P} are being approximated, the surrogate solutions $\tilde{\mathbf{x}}$ obtained by solving~\eqref{alg:relaxation} may not be feasible for~\eqref{sec:surr_problem:P}. Indeed, none of the surrogate solutions found in Table~\ref{tab:WDN_hanoi_surrogate} is feasible for~\eqref{sec:surr_problem:P}. Thus, MISSOC proceeds to the post processing phase in Algorithm~\ref{alg:postprocessing}, which is illustrated in the next subsection.

{
\renewcommand{\arraystretch}{1.2}
\begin{table}[h!]
 \centering
\begin{tabular}{@{}crrr@{}}
\toprule

       \textbf{Metric} & \multicolumn{1}{c}{\textbf{Gurobi}} & \multicolumn{1}{c}{\textbf{BARON}} & \multicolumn{1}{c}{\textbf{SCIP}}\\ \toprule
         $g_0(\tilde{\mathbf{x}})$   & 5750110.400&  5747152.900  & n/a  \\
         gap(\%) & 0&  0  &  n/a      \\
         time(s)  & 13.773  &  153.635 &  600    \\ 
         \bottomrule
\end{tabular}
\caption{Results for the relaxed surrogate problem~\eqref{alg:relaxation} of the Hanoi instance.}
\label{tab:WDN_hanoi_surrogate}
\end{table}
}

\subsection{Post processing}
\label{sec:demonstrative:postprocessing}

In this subsection, we illustrate the post processing step of MISSOC, summarized in Algorithm~\ref{alg:postprocessing}. Starting from the surrogate solution $\tilde{\mathbf{x}}$ obtained as described in Section~\ref{sec:demonstrative:surrogate}, the goal is to recover feasibility with respect to the original problem~\eqref{sec:surr_problem:P} and further improve the solution. Since the Water Distribution Network problem is a design problem, the main goal is to obtain the best possible solution rather than a fast feasible one. Therefore, in this example, we  perform the optional third subphase of the post processing step, which focuses on improving the objective value.

To simplify the presentation, we focus on the surrogate solution $\tilde{\mathbf{x}}$ obtained with Gurobi in Section~\ref{sec:demonstrative:surrogate} and illustrate the post processing phase only from this starting point. We also keep Gurobi as the solver within the post processing phase. To make a fair comparison with the solution of the original problem, the overall time limit for MISSOC is set to 600 seconds. Since solving the relaxed surrogate problem ~\eqref{alg:relaxation}  with Gurobi takes 13.773 seconds, the remaining time available for the post processing step is 586.227 seconds.

As explained in Algorithm~\ref{alg:postprocessing}, the weights $\omega_0$ and $\omega_u$ in the first subphase are chosen by the user. In the Water Distribution Network problem, the objective function $g_0$ (total cost of the water network) is completely determined by the binary variables that decide which diameter is assigned to each pipe. Therefore, when these binary variables are fixed, the value of $g_0$ is constant. For this reason, in the first subphase we set $\omega_0=0$ and $\omega_u=1$, so that the subproblem focuses on reducing the feasibility violation measured by the slack variable $u$. Moreover, in the optional third subphase, we set the relative change tolerance to $\tau=0.01$ and the patience parameter to $\eta=2$.

Table~\ref{tab:WDN_hanoi_postprocessing} reports the results obtained at each iteration $\kappa$ during the post processing step, namely Algorithm~\ref{alg:postprocessing}, indicating the corresponding subphase. For each iteration, we report the objective value of the solution found, $g_0(\mathbf{x}^{(\kappa)})$, the computational time in seconds and the corresponding slack value $u^{(k)}$. The latter is used to identify when feasibility has been recovered, i.e., when $u^{(k)}\leq 10^{-6}$. We also report the accumulated computational time in MISSOC (in seconds), denoted by $T$, computed as the time spent to obtain $\tilde{\mathbf{x}}$ in Section~\ref{sec:demonstrative:surrogate} with Gurobi, 13.773 seconds, plus the accumulated time spent in the post processing step. This allows us to track the quality and feasibility of the solutions returned by MISSOC as the computational time increases. 

Here, $\kappa$ represents the number of pipes for which the selected diameter is allowed to change. In the extended formulation used here, the diameter selection for each pipe is modeled with one binary variable for each pipe and each available diameter. Exactly only one of these variables must be equal to one for each pipe. Thus, changing the selected diameter of one pipe necessarily requires two binary switches in the extended formulation: one diameter option is deselected and another one is selected for that same pipe. Therefore, the corresponding local branching constraint in the extended binary variables needs to allow up to $2\kappa$ binary switches, so that the diameter selection can change in up to $\kappa$ pipes.

As shown in Table~\ref{tab:WDN_hanoi_postprocessing}, feasibility is recovered when $\kappa=3$, meaning that the selected diameter is changed in three pipes with respect to the initial solution $\tilde{\mathbf{x}}$. These changes involve less than $10\%$ of the pipes in the network. In terms of the extended formulation, this corresponds to changing $2\kappa=6$ binary variables, i.e., approximately $3\%$ of the 204 binary variables. This first feasible solution is obtained after 21.823 seconds of total MISSOC time and has an objective value of 6442135.300. It is worth comparing this result with the first feasible solution found when solving the original problem~\eqref{sec:surr_problem:P} directly with Gurobi. In that case, the first feasible solution has a worst objective value of 6595608.200 and is found after approximately 45 seconds. Thus, MISSOC obtains its first feasible solution faster and with a better objective value than the first feasible solution found by directly solving~\eqref{sec:surr_problem:P}.

After feasibility is recovered, MISSOC proceeds to the optional third subphase to further improve the objective value while maintaining feasibility. The objective value stabilizes at $\kappa=6$ and the corresponding solution is returned as the best solution found.  This means that changing the diameter selection in 6 of the 34 pipes is enough to improve the solution until stabilization. This corresponds to  $2\kappa=12$ binary switches in the extended formulation, i.e., less than $6\%$ of the total 204 binary variables.

The best solution found by MISSOC has objective value of $g_0(\mathbf{x}^*)=6109620.899$ and is obtained after 248.644 seconds in total, including both the time required to solve~\eqref{alg:relaxation} with Gurobi and the time spent in the post processing step. Observe that this objective value coincides with the optimal value obtained by directly solving ~\eqref{sec:surr_problem:P} (see Table~\ref{tab:WDN_hanoi_original}). Therefore,  MISSOC reaches the same optimal value in less time than directly solving the extended formulation, which takes at least 386.370 seconds with Gurobi. However, we recall that MISSOC cannot provide optimality guarantees for the original problem.

{
\renewcommand{\arraystretch}{1.2}
\begin{table}[h!]
 \centering
\begin{tabular}{@{}lcrrrr@{}}
\toprule
\textbf{subphase} &
\multicolumn{1}{c}{$\boldsymbol{\kappa}$} &
\multicolumn{1}{c}{$\mathbf{g_0(\mathbf{x}^{(\kappa)})}$} &
\multicolumn{1}{c}{\textbf{time (s)}} & \multicolumn{1}{c}{\textbf{slack value} $u$} & 
\multicolumn{1}{c}{$\mathbf{T}$ \textbf{(s)}} 
 \\\toprule
         First   & 0&  5750110.400   &  0.120&0.120 &13.894  \\
         Second   & 1& 5870072.300   & 0.870 &0.051 & 14.763 \\
         Second   & 2&  6285121.300  &  3.310  &0.009& 18.073 \\
         Second   & 3&   6442135.300 &  3.750   & $\leq 10^{-6}$ &   21.823  \\
         Third    & 4&     6272729.800  &  29.110 &$\leq 10^{-6}$ & 50.933      \\
         Third    & 5&      6167726.300 &  61.480 & $\leq 10^{-6}$ & 112.413     \\
          Third   & 6&  6109620.899  &  136.231 &$\leq 10^{-6}$ & 248.644       \\
         \bottomrule
\end{tabular}
\caption{Post processing results for the Hanoi instance starting from $\tilde{\mathbf{x}}$ obtained with Gurobi as described in Section~\ref{sec:demonstrative:surrogate}.}
\label{tab:WDN_hanoi_postprocessing}
\end{table}
}

This case study shows that MISSOC can find high-quality feasible solutions for real-world MINLPs with complicating constraints and integer variables. It illustrates how the post processing step can recover feasibility and improve the objective value within a limited computational time, starting from the solution of the relaxed surrogate problem \eqref{alg:relaxation}. Moreover, the first feasible solution found by MISSOC is obtained faster and has a better objective value than the first feasible solution obtained by directly solving~\eqref{sec:surr_problem:P}. This can be relevant in practical settings where obtaining a good feasible solution quickly may be more important than proving optimality.

\section{Computational results}
\label{sec:comput_results}

In this section, we further evaluate the proposed MISSOC algorithm through a set of computational experiments. We consider both benchmark instances from the MINLPlib library \citep{bussieck2003minlplib} and the real-world case study of the Hydro Unit Commitment problem \citep{borghetti2015optimal,Taktak2017}.

We emphasize that MISSOC is an approach aimed at addressing challenging MINLPs in which off-the-shelf solvers struggle to find feasible good-quality solutions. Note that benchmark libraries such as MINLPlib are widely used in the development, enhancement and evaluation of global optimization solvers. As a result, their instances are well studied and modern solvers often perform well on them.  For this reason, MINLPlib instances are used in this work as a controlled setting to evaluate MISSOC. In contrast, in the real-world case study considered in this section, solvers face significantly greater challenges. Therefore, MISSOC can more effectively demonstrate its full potential for solving complex MINLPs in realistic settings.

This section is organized as follows. We first describe the experimental setup in Section~\ref{sec:exp:setup}. Then, in Section~\ref{sec:exp:MINLPlib}, we present the results obtained on the MINLPlib benchmark instances. Finally, in Section~\ref{sec:exp:HUC}, we evaluate MISSOC in the real-world Hydro Unit Commitment problem.

\subsection{Set up}
\label{sec:exp:setup}

We describe the experimental setup used in all computational experiments to ensure consistency and reproducibility across the different instances.

For all tested instances, the size of the training set used to approximate each function is determined based on the complexity of the smooth additive regression model employed. The complexity of a regression model is determined by the number of its coefficients. In our case, it is the sum of degrees $d_j$ and the number of intervals $k_j$ across all $p$ covariates, plus one accounting for the intercept $\alpha$. To ensure a reliable fit, we generate 15 samples per model coefficient. Therefore, the total number of observations $n$ in the training dataset is $15 \left( 1+\sum_{j=1}^p(d_j + k_j)\right)$.

For fitting a smooth additive regression model using the previously generated training dataset, MISSOC assumes that the degree $d_j$ and the number of intervals $k_j$ of each variable $x_j$ are given as inputs (see Algorithm~\ref{alg:missoc}).  These parameters can be selected in different ways. One option could be to perform a grid search over a range of candidate degrees and number of intervals, selecting the best combination according to an out-of-sample regression metric. However, this process can be time consuming and is often not strictly necessary. Since the goal is to build a surrogate problem that approximates the original one to find a good solution, a good enough approximation of the complicating function is often practical, without needing to be the best one over a candidate grid. For this reason, in practice, we set the degree and number of intervals based on prior experience, using values that have performed well in preliminary computational studies. In our experience, using 10 intervals and degree 3 for all variables provides a good balance between approximation accuracy and solution quality, reducing the computational costs of the approximating phase. Thus, these values are used throughout this work. In addition, the knots defining the intervals are placed equidistantly.

If there is expert knowledge about the shape of the complicating function to be approximated, we fit a shape-constrained model by incorporating the corresponding constraints  as described in Section~\ref{sec:shape_constrained_regression}. Otherwise, we fit an unconstrained model by solving problem~\eqref{eq:preliminaries:PLSE_additive_matrix} described in Section~\ref{sec:preliminaries:smooth_additive_regresion}. In both cases, the resulting optimization problem is solved with MOSEK using the Python Fusion API \citep{mosek}. The implementation was carried out in Python v.\,3.12.8.  For reproducibility, an open source repository containing the implementation for fitting these smooth additive regression models (with or without constraints) is available at: \url{https://github.com/RSpencerT/MISSOC/tree/main/shape_constrained_smooth_additive_regression_model}.  All approximations performed in this work can be reproduced using the provided datasets and scripts. 

The resulting surrogate problem~\eqref{sec:surr_problem:tildeP} is solved using different solvers for comparison purposes. First, we use the customized SC-MINLP algorithm \citep{d2012algorithmic,d2019strengthening} (using CPLEX v.\,22.1.1), which leverages the additive and separable nature of the approximating function formulation. Besides, we also use Gurobi v.\,12.0.1 \citep{gurobi}, BARON v.\,24.5.8 \citep{sahinidis1996baron} and COUENNE v.\,0.5.8 \citep{belotti2009branching}. For the post processing phase, the choice of solver depends on the subphase considered. In the first subphase, the problem is continuous and is solved with the local solver IPOPT v.\,3.14.14 \citep{wachter2006implementation}. If the subsequent local branching subphases are performed, an MINLP solver such as Gurobi is required.

We evaluate MISSOC by comparing its heuristic solution $\mathbf{x}^*$ (see Algorithm~\ref{alg:missoc}) with the solutions obtained by directly solving the original problems~\eqref{sec:surr_problem:P} with the same global solvers. We record the computation time (in seconds) and the objective value of the achieved solution. Note that the computational time in MISSOC includes the time required for solving the surrogate problem~\eqref{sec:surr_problem:tildeP} plus the time spent in the post processing step. To compare the performance of the solvers in solving~\eqref{sec:surr_problem:tildeP} (i.e., before the post processing step), we also record the optimality gap. We compute it using the formula implemented in Gurobi v.\,12.0.1:

\begin{equation}
    \text{optimality gap} = 100 \cdot \frac{|\text{upper bound} - \text{lower bound}|}{|\text{upper bound}|},
    \label{eq:optimality_gap}
\end{equation}

\noindent where the upper bound and lower bound are those reported by the corresponding solver. The optimality gap obtained when solving the surrogate ~\eqref{sec:surr_problem:tildeP} must not be compared to that when solving the original problem~\eqref{sec:surr_problem:P}, since the final solution $\mathbf{x}^*$ provided by MISSOC, i.e., after the post processing step, is heuristic and, thus, we cannot provide optimality guarantees.

To solve both the original and surrogate problems, all solvers have been run
with a time limit of 600 seconds and the rest of the parameters are set to default. The experiments have been conducted on a machine with an Intel(R) Xeon(R) CPU E5-2620 v4 @ 2.10GHz processor and 64GB of RAM.

\subsection{Benchmark MINLPLib instances}
\label{sec:exp:MINLPlib}

Three instances from the MINLPlib library are considered to benchmark MISSOC: \textit{ex6\_2\_13}, \textit{ex6\_2\_5} and \textit{ex6\_2\_7}. These have been selected from the test set in \cite{bertsimas2025global} with the following criteria. First we selected  the problems where the complicating function  appears in the objective. Among them, we consider only those that are challenging for Gurobi v.\,12.0.1, i.e., those that it cannot solve in 600 seconds. Let us highlight that \cite{bertsimas2025global} used Gurobi v.\,8 in their experiments. Since Gurobi v.\,12 includes better capabilities for solving global optimization problems, many of their benchmark instances have been excluded from our study because they are now easily solved. 

MINLPlib reports the best feasible solution known for each of the selected instances: \textit{ex6\_2\_13} ($-0.216$)\footnote{\url{https://www.minlplib.org/ex6_2_13.html}}, \textit{ex6\_2\_5} ($-70.752$)\footnote{\url{https://www.minlplib.org/ex6_2_5.html}} and \textit{ex6\_2\_7} ($-0.161$)\footnote{\url{https://www.minlplib.org/ex6_2_7.html}}. However, they are not marked as solved to global optimality. The characteristics of these selected instances are summarized in Table~\ref{tab:characterisics_MINLPlib}. The first column indicates the instance name. The second column indicates whether the information corresponds to the original problem~\eqref{sec:surr_problem:P}, as formulated in MINLPlib, or to the surrogate problem~\eqref{sec:surr_problem:tildeP} that is built later by MISSOC.  The third, fourth and fifth columns report the number of constraints, binary variables and continuous variables, respectively. It can be observed that the selected problems are all continuous.

{
\renewcommand{\arraystretch}{1.35}
\begin{table}[h!]
 \centering
 \small
\begin{tabular}{ccccc}
\toprule
\textbf{Instance} & \textbf{Problem} & \textbf{\# constraints}  & \textbf{\# binary vars} & \textbf{\# continuous vars} \\\toprule
\multirow{2}{*}{\textit{ex6\_2\_13}} 
  &\eqref{sec:surr_problem:P}   & 3 & 0 & 6 \\
  &\eqref{sec:surr_problem:tildeP}  & 201  & 60  &  132  \\ \midrule
\multirow{2}{*}{\textit{ex6\_2\_5}}
  &\eqref{sec:surr_problem:P}  & 3 & 0 & 9 \\
  &\eqref{sec:surr_problem:tildeP} &  300 &  90 &   198 \\ \midrule
\multirow{2}{*}{\textit{ex6\_2\_7}}
  &\eqref{sec:surr_problem:P}  & 3 & 0 & 9 \\
  &\eqref{sec:surr_problem:tildeP}& 300  &  90  &  198  \\ \bottomrule
\end{tabular}
\caption{Number of constraints, binary variables and continuous variables in the original~\eqref{sec:surr_problem:P} and surrogate~\eqref{sec:surr_problem:tildeP}  problems of the MINLPlib instances.}
\label{tab:characterisics_MINLPlib}
\end{table}
}

The objective function in each of the three benchmark instances includes both linear and nonlinear terms. To build the surrogate problems, only the nonlinear ones are approximated while the linear terms remain unchanged. Since no expert knowledge is available, the nonlinear terms are approximated with a smooth additive regression model without shape constraints through problem~\eqref{eq:preliminaries:PLSE_additive_matrix}.  The fitting process is fast, requiring a CPU time in the order of $10^{-2}$ seconds per instance.

The characteristics of the surrogate problems of the MINLPlib instances are detailed in Table~\ref{tab:characterisics_MINLPlib}.  As explained in Section~\ref{sec:surrog_formulation},  each knot in each covariate involved in the regression model  to approximate the objective function introduces an additional binary variable in the formulation of the surrogate problem, as well as other continuous variables. A set of constraints also needs to be incorporated. As a consequence, the size of the surrogate problem is always larger than the original. 

Since the original MINLPlib instances considered here are continuous and the approximated function appears in the objective, only the first subphase of the post processing step in MISSOC can be applied. No slack variable need to be introduced in problem~\eqref{alg:ls_problem} and the post processing problem is solved with $\omega_0=1$.

The results for the three MINLPlib benchmark instances are presented in Table~\ref{tab:results_MINLPlib}. The first column shows the instance name. The second column indicates the method used to solve each instance, distinguishing between directly solving~\eqref{sec:surr_problem:P} with a MINLP solver and applying MISSOC. The third column lists the different evaluation metrics that are reported, which are detailed next. First, the value of the original objective function evaluated at the achieved solution is reported. When directly solving the original problems, this is $g_0(\mathbf{x}^{(P)})$, being $\mathbf{x}^{(P)}$ the corresponding solution found.  When applying MISSOC, it is $g_0(\mathbf{x}^*)$, i.e., evaluated at the heuristic solution.
Next, the computational time in seconds is reported (time(s)), followed by the optimality gap in percent (gap(\%)) computed as in~\eqref{eq:optimality_gap}. The remaining columns provide the values of these metrics for each solver: SC-MINLP, Gurobi, BARON and COUENNE. 
Recall that the SC-MINLP algorithm is tailored for solving optimization problems where the nonconvexities appear only in additive and univariate structures. Since the original problems of the MINLPlib instances do not have such a structure, their corresponding entries are marked with ``–''.

{
\renewcommand{\arraystretch}{1.35}
\begin{table}[h!]
 \fontsize{8.7pt}{11pt}\selectfont
 \centering
\begin{tabular}{@{}cccrrrr@{}}
\toprule
        \textbf{Instance} & \textbf{Method} & \textbf{Metric} &\textbf{SC-MINLP} & \textbf{Gurobi} & \textbf{BARON} &  \textbf{COUENNE}\\ \toprule
       \multirow{7}{*}{\textit{ex6\_2\_13}} & \multirow{3}{*}{MINLP solver } & $g_0(\mathbf{x}^{(P)})$ & -&  $-0.216$&  $-0.216$&   $-0.216 $\\
        & & time(s) & -& 600.205& 603.305& 607.255 \\
        & & gap(\%) & -& 18.603 & 5.243 & 407.699 \\ 
        \cmidrule(l){2-7}
      & \multirow{3}{*}{\makecell[c]{MISSOC }} & $g_0(\mathbf{x}^*)$ & $-0.216$& $-0.216$& $-0.216$& $-0.216$ \\
       &  & time(s) & 0.028& 0.823& 603.152&  3.441\\
    & & gap(\%) & 0 & 53.008 \textsuperscript{\dag} &5.027& 0 \\  \midrule
       
    \multirow{7}{*}{\textit{ex6\_2\_5}} & \multirow{3}{*}{ MINLP solver} & $g_0(\mathbf{x}^{(P)})$ & -& $-70.599$ & $-70.752$&  $-70.752$\\
        & & time(s) & -&600.278 & 603.024&  607.002\\
         & & gap(\%) & -& 345.257& 203.848& 1386.825   \\ 
        \cmidrule(l){2-7}
      & \multirow{3}{*}{\makecell[c]{MISSOC}} & $g_0(\mathbf{x}^*)$ & $-70.558$& $-70.558$& $-70.558$&  -70.552\\
       &  & time(s) & 0.127 & 0.456& 14.620& 63.842 \\
       & & gap(\%) & 0.005 & 0.001& 0& 0 \\ \midrule

    \multirow{7}{*}{\textit{ex6\_2\_7}} & \multirow{3}{*}{ MINLP solver } & $g_0(\mathbf{x}^{(P)})$ &- &$-0.161$ &$-0.161$ &$-0.161$  \\
        & & time(s) & - & 600.633& 602.975& 607.277  \\
        & & gap(\%) & -& 2400.913& 1359.993& 10913.441  \\ 
        \cmidrule(l){2-7}
      & \multirow{3}{*}{\makecell[c]{MISSOC}} & $g_0(\mathbf{x}^*)$ &$-0.161$ & n/a & $-0.161$&  $-0.161$\\
       &  & time(s) &0.115 & 600.312& 602.997& 607.772 \\
       & & gap(\%) & 0.016& n/a & 14.321& 17.842 \\ \bottomrule     
\end{tabular}
\caption{Results obtained by directly solving the original MINLPlib instances and by solving them with MISSOC.}
\label{tab:results_MINLPlib}

{\raggedright \textsuperscript{\dag} Gurobi terminated with the message: Sub-optimal termination (unable to solve some node relaxations). \par}
\end{table}
}

To assess the quality of the solutions, we first compare the objective values obtained  using MISSOC with the above-mentioned best known feasible solutions reported by MINLPlib for each case. The three instances are minimization problems, so lower objective values indicate better solutions. As seen in Table~\ref{tab:results_MINLPlib}, the objective values returned by MISSOC  are either the same or very close to the reference values in all cases.

Regarding computational times, we observe that solving the original problems directly with Gurobi, BARON and COUENNE always requires more than 600 seconds. In contrast, addressing the problems with MISSOC results in lower computational times in many cases. However, the solvers in the original problem~\eqref{sec:surr_problem:P} often find good feasible solutions quickly and spend the remaining time trying to close the optimality gap. To illustrate this, we compare the solutions obtained by the solvers for~\eqref{sec:surr_problem:P} at the time spent by MISSOC with SC-MINLP to return its solution. For \textit{ex6\_2\_13}, after 0.028 seconds, all solvers obtain the $-0.216$ objective value.  For \textit{ex6\_2\_5}, after 0.127 seconds, BARON and COUENNE already reach the best known value of $-70.752$ and Gurobi is close. For \textit{ex6\_2\_7}, after 0.115 seconds, all three solvers obtain the $-0.161$ objective value.

Comparing the times spent by MISSOC to return the heuristic solution $\mathbf{x}^*$ using the different solvers, SC-MINLP stands out in comparison to the other solvers. When using SC-MINLP to solve the surrogate problem, MISSOC always take   less than 0.15 seconds to provide the heuristic solution. With the other solvers, the computational times vary between instances. For example, using Gurobi in MISSOC is fast on \textit{ex6\_2\_13} and 
\textit{ex6\_2\_5} (under one second), but it reaches the time limit without finding a feasible solution for \textit{ex6\_2\_7} (indicated with ``n/a'' in the table). One can note that even in its best cases, Gurobi remains slower than SC-MINLP. When using BARON and COUENNE  to solve the surrogate problems, MISSOC takes moderate time in some instances, but can require 600 seconds on other instances. In the particular case of the surrogate \textit{ex6\_2\_7} instance, only SC-MINLP solves it before the time limit.

Overall, the results obtained on these benchmark instances confirm that MISSOC can provide high-quality solutions close to the best known feasible ones in short computational times, especially when the surrogate problem is solved with SC-MINLP. Indeed, SC-MINLP yields the most consistent results across the three instances considered. It achieves solutions very close to the best known feasible ones with low computational effort and zero or near-zero optimality gaps in the surrogate problems.

\subsection{Real case study: Hydro Unit Commitment problem}
\label{sec:exp:HUC}

We now evaluate MISSOC on the well-known real-world Hydro Unit Commitment problem \citep{borghetti2015optimal,Taktak2017}. It involves determining the optimal schedule for operating hydroelectric units over a given time horizon to maximize revenue subject to technical and operational constraints.  The formulation and a feasible solution with an objective value of 14533.1 are reported in \cite{borghetti2015optimal}.

This problem is particularly relevant for assessing MISSOC, since it is a realistic and challenging nonconvex problem with integer decision variables that  commercial solvers struggle to solve. In particular, the complicating nonlinearities arise in the power generation function $P(v,w)$, which appears in the objective function and depends on reservoir volume $v$ and water flow $w$. In this case, expert knowledge indicates that the power generation output lies between 2.595 and 24.089. These values correspond to the power generation function evaluated at the minimum and maximum values of $v$ and $w$, respectively. 
To ensure the approximation aligns with the known behavior, these values are imposed as lower and upper bounds in the smooth additive regression model that approximates the true function. Then, we apply the methodology for shape-constrained estimation described in Section~\ref{sec:shape_constrained_regression} and estimate the model  through problem~\eqref{eq:shape_constrained_additive:shape_constrained_problem}.  This illustrates the knowledge-driven aspect of MISSOC, since available shape information is incorporated into the data-driven surrogate.

\begin{figure}[h!]
    \centering
    \begin{subfigure}{0.32\textwidth}
        \includegraphics[width=\linewidth]{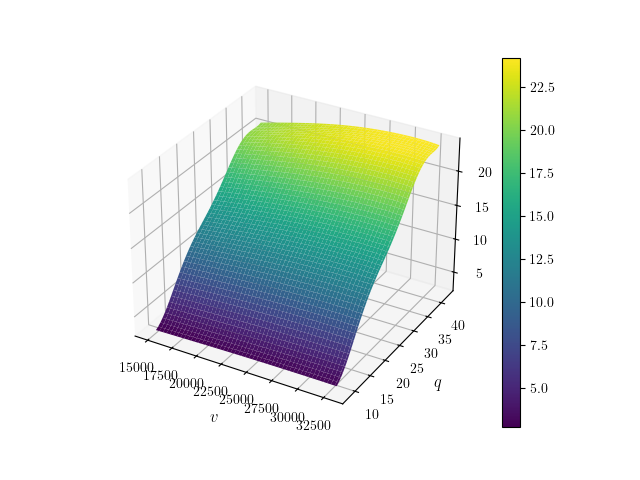}
        \caption{Real surface.}
    \end{subfigure}
    \hfill
    \begin{subfigure}{0.32\textwidth}
        \includegraphics[width=\linewidth]{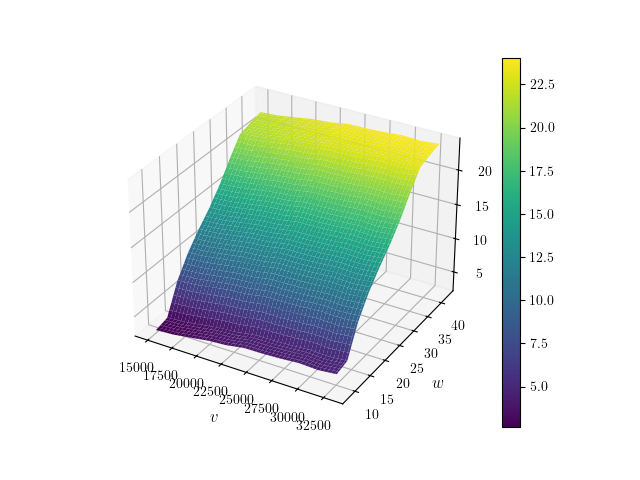}
        \caption{$d_j = 1, \ j=1,2$.}
    \end{subfigure}
    \hfill
    \begin{subfigure}{0.32\textwidth}
        \includegraphics[width=\linewidth]{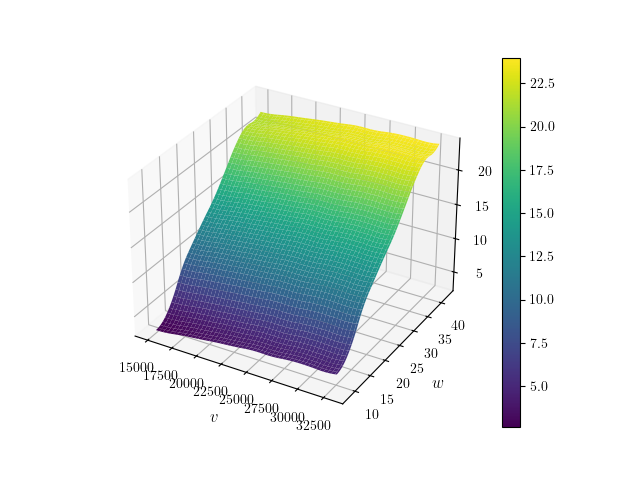}
        \caption{$d_j = 2, \ j=1,2$.}
    \end{subfigure}
      \hfill
    \begin{subfigure}{0.32\textwidth}
    \includegraphics[width=\linewidth]{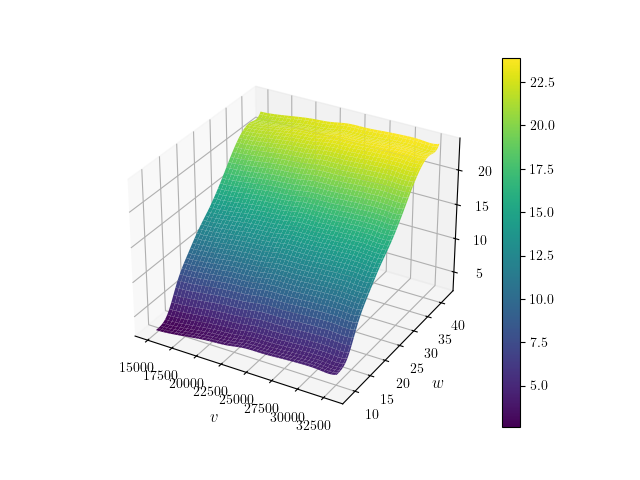}
    \caption{$d_j = 3, \ j=1,2$.}
\end{subfigure}
  \hfill
    \begin{subfigure}{0.32\textwidth}
    \includegraphics[width=\linewidth]{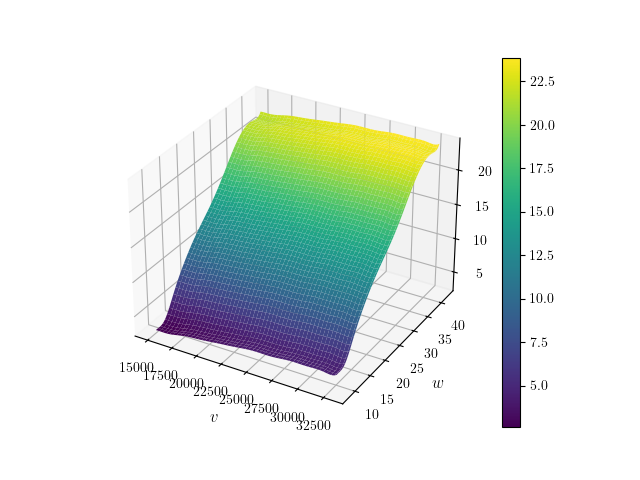}
    \caption{$d_j = 4, \ j=1,2$.}
\end{subfigure}
  \hfill
 \begin{subfigure}{0.32\textwidth}
    \includegraphics[width=\linewidth]{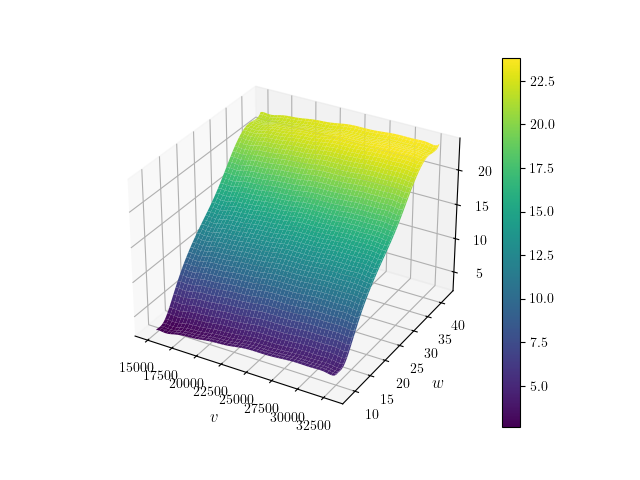}
    \caption{$d_j = 5, \ j=1,2$.}
    \end{subfigure}
  \hfill
 \begin{subfigure}{0.32\textwidth}
    \includegraphics[width=\linewidth]{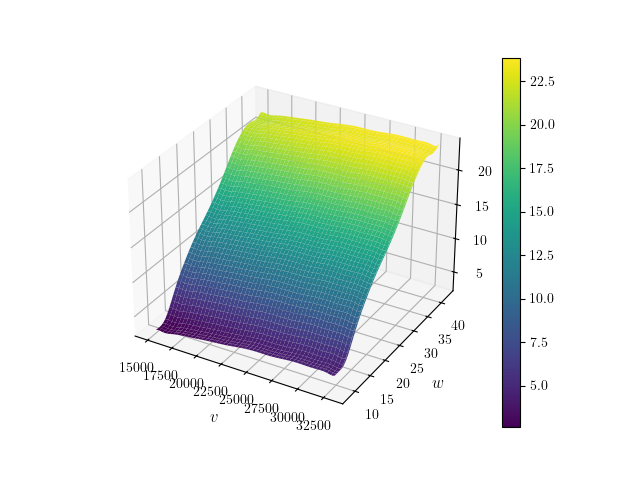}
    \caption{$d_j = 6, \ j=1,2$.}
\end{subfigure}
 \begin{subfigure}{0.32\textwidth}
    \includegraphics[width=\linewidth]{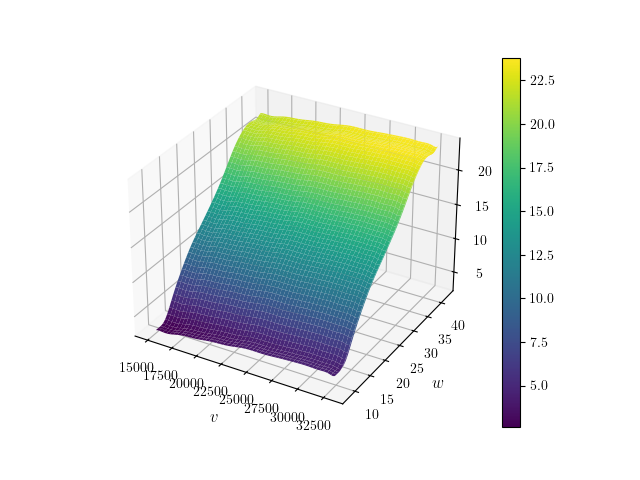}
    \caption{$d_j = 7, \ j=1,2$.}
\end{subfigure}
\caption{Power generation surface \( P(v, w) \) in the Hydro Unit Commitment problem and its bounded approximations using degrees 1 to 7 in both covariates \( v \) and \( w \), with 10 intervals per variable.}
    \label{fig:hydro}
\end{figure}

\begin{figure}[t]
    \centering

    \begin{subfigure}[t]{0.49\textwidth}
        \centering
        \includegraphics[width=\textwidth]{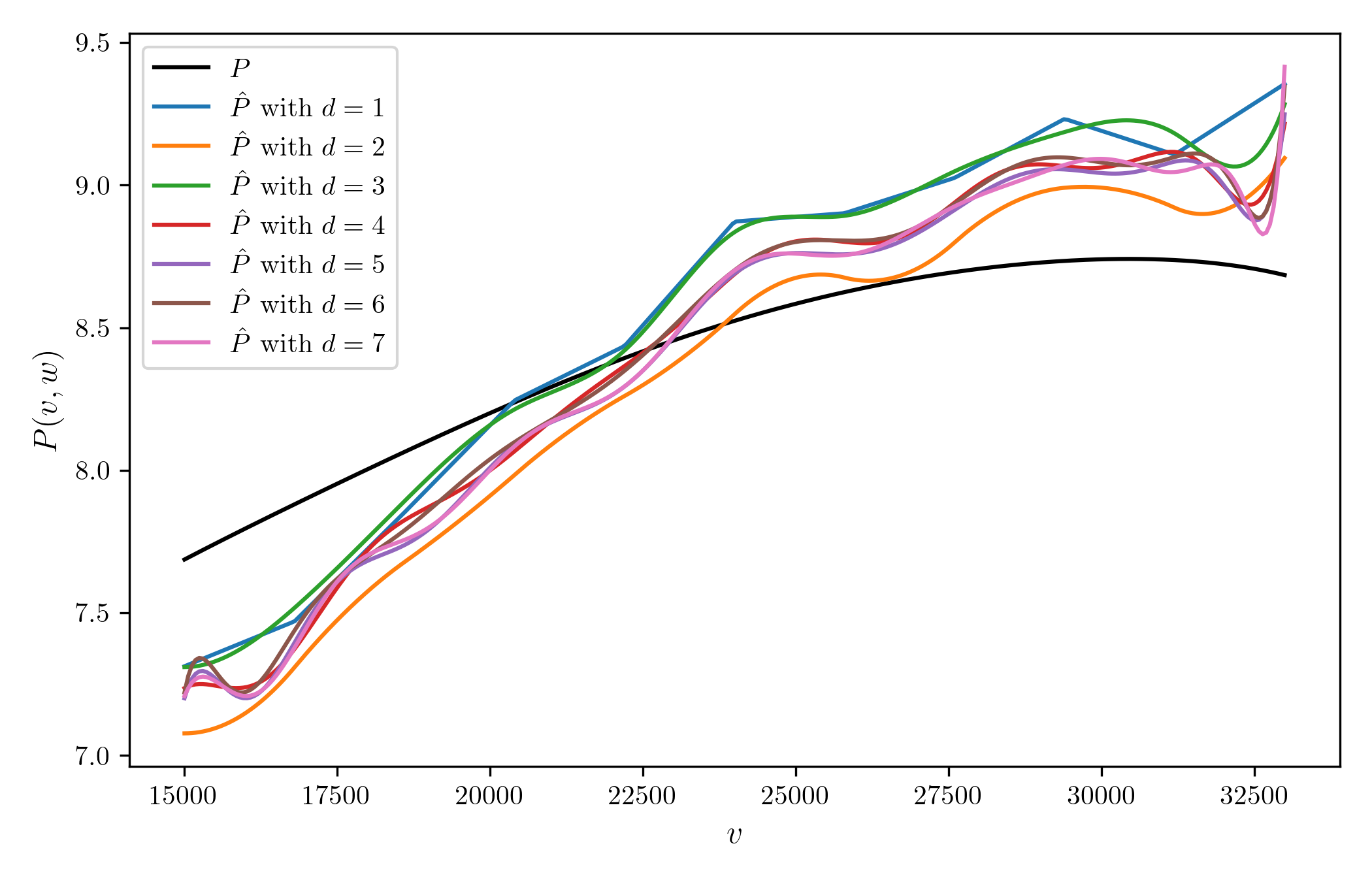}
        \caption{$w$ fixed at the 25th percentile (15.42), varying $v$.}
    \end{subfigure}
    \hfill
    \begin{subfigure}[t]{0.49\textwidth}
        \centering
        \includegraphics[width=\textwidth]{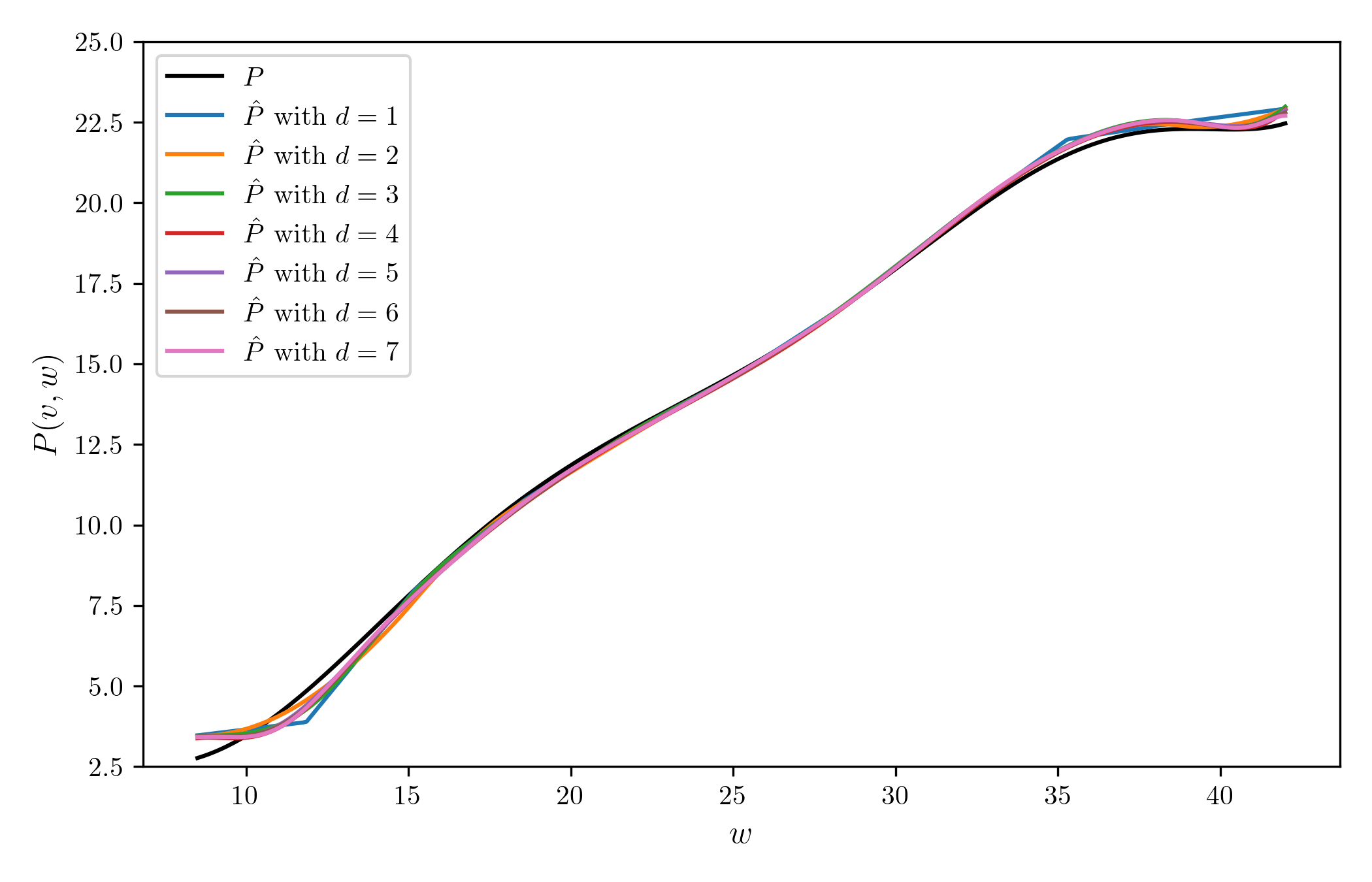}
        \caption{$v$ fixed at the 25th percentile (20139.58), varying $w$.}
    \end{subfigure}

\caption{One dimensional slices of the true function $P(v,w)$ in the Hydro Unit Commitment problem and of its fitted approximating functions for degrees 1 to 7.}
    \label{fig:huc_slices_all}
\end{figure}

To illustrate the effect of the choice of the degree in the solution 
provided by MISSOC, we fit models using degrees 1 to 7. In all cases, the computational time for fitting the models is below 1.5 seconds, with lower-degree models requiring slightly less time.
Figure~\ref{fig:hydro} includes the visualization of the real power generation function $P(v,w)$ together with its seven different approximations. As a complement and to deeply explore the effect of the degree, we have obtained one dimensional slices of the true function $P(v,w)$ and of its fitted approximating functions using degrees 1 to 7. For this, one coordinate is fixed at a specific value and the other is allowed to vary. We have analyzed slices corresponding to different fixed values of the two coordinates and we have observed that the effect of the degree is more clearly visible when fixing $w$ and varying $v$ than in the opposite direction. As an illustrative example, Figure~\ref{fig:huc_slices_all} shows representative slices obtained by fixing $w$ and $v$ at their respective 25th percentiles.

Next, a surrogate problem is built for each of the seven approximations of the objective function.  The different steps for building them, as described in Section~\ref{sec:surrog_surrogate-problem}, and solving them using the SC-MINLP algorithm are detailed in the following Jupyter Notebook publicly available at GitHub: \\ \url{https://github.com/RSpencerT/MISSOC/blob/main/surrogate_optimization/Splines_HUC.ipynb}.

Table~\ref{tab:characterisics_HUC} summarizes the characteristics of the original and surrogate problems in the Hydro Unit Commitment instance, following the same structure as in Table~\ref{tab:characterisics_MINLPlib}. Since the number of knots remains fixed across all approximations and only the degree changes, the seven surrogate problems corresponding to degrees 1 to 7 have the same sizes.

{
\renewcommand{\arraystretch}{1.35}
\begin{table}[h!]
 \centering
\begin{tabular}{@{}cccc@{}}
\toprule
\textbf{Problem}                             & \multicolumn{1}{c}{\textbf{\# constraints}}  & \textbf{\# binary vars} & \textbf{\# continuous vars} \\ \toprule
\multicolumn{1}{c}{\eqref{sec:surr_problem:P}}        &    169    &         72    &    48     \\
\multicolumn{1}{c}{\eqref{sec:surr_problem:tildeP}}    &         1201   &       312 &      552         \\ \bottomrule
\end{tabular}
\caption{Number of constraints, binary variables and continuous variables in the original~\eqref{sec:surr_problem:P}  and ~\eqref{sec:surr_problem:tildeP} problem of the Hydro Unit Commitment instance.}
\label{tab:characterisics_HUC}
\end{table}
}

The Hydro Unit Commitment problem is an operational problem, where obtaining good feasible solutions quickly is prioritized over spending additional time on further objective improvement. Therefore, to keep the setting realistic, we do not perform the optional third subphase of the post processing step, designed for objective improvement through local branching. In addition, since the objective function is approximated but not the constraints, the feasibility recovery through the second subphase is not required. As a result, only the first subphase of the post processing step is applied. Since no constraint is approximated, no slack variable is introduced and problem~\eqref{alg:ls_problem} is solved with $\omega_0=1$.

The results for the Hydro Unit Commitment problem are reported in Table~\ref{tab:results_HUC}, which has a similar structure to Table~\ref{tab:results_MINLPlib}. Since we vary the degree of the $B-$spline basis for the approximation of the objective function, when using MISSOC, we also report $\tilde{g}_0( \mathbf{\tilde{x}})$, i.e., the value of the approximated objective function evaluated at the surrogate solution $\mathbf{\tilde{x}}$, which is already feasible to~\eqref{sec:surr_problem:P} in this case. This allows us to evaluate how the approximation behaves as the degree increases. Moreover, since the Hydro Unit Commitment problem originally contains integer and continuous variables, we also test here the BONMIN v.\,1.8.9 solver \citep{bonami2008algorithmic}. However, note that BONMIN cannot provide a reliable dual bound when solving nonconvex problems. As a result, the optimality gaps are not reported and it is indicated with ``–'' in the table. Recall that the SC-MINLP algorithm cannot be used to solve the original problem because it does not have a separable and additive structure. Then, the corresponding entries are also filled with ``–''.

{
\renewcommand{\arraystretch}{1.2}
\begin{table}[h!]
 \centering
  \fontsize{8.7pt}{11pt}\selectfont
\begin{tabular}{@{}ccrrrrr@{}}
\toprule

         \textbf{Method} & \textbf{Metric} &\textbf{SC-MINLP} & \textbf{Gurobi} & \textbf{BARON} & \textbf{COUENNE}& \textbf{BONMIN}\\ \toprule

         \multirow{3}{*}{MINLP solver} & $g_0(\mathbf{x}^{(P)})$ & - & \multirow{3}{*}{infeas.} & 10350.291& 6600.633 & 7814.439\\
         & time(s) & - & & 603.103& 606.776 &  2.092\\
         & gap(\%) & -& & 72.407& 71.886& -\\ \midrule
        \multirow{4}{*}{\makecell[c]{ MISSOC \\ with  degree 1}} &$g_0(\mathbf{x}^*)$ &14519.363
& 14519.363 & 14490.601 & 14519.363&14519.363 \\
        & $\tilde{g}_0( \mathbf{\tilde{x}})$ & 14790.850& 14790.850& 14790.850 & 14790.850 &14790.850 \\
         & time(s) & 0.210& 0.195& 0.792& 602.758& 606.969\\
         & gap(\%) & 0& 0& 0& 17.536& -\\\midrule
        \multirow{4}{*}{\makecell[c]{ MISSOC \\ with  degree 2}}  & $g_0(\mathbf{x}^*)$ & 14519.363 &14519.363 & 14519.363& 14127.634& 14392.289\\
        & $\tilde{g}_0( \mathbf{\tilde{x}})$ & 14617.749
& 14617.721 &14617.749& 14153.724& 14456.557\\
         & time(s) & 0.557& 2.018& 149.820
&602.619 & 607.567\\
         & gap(\%) & 0& 0.010& 0& 39.234& -\\\midrule
        \multirow{4}{*}{\makecell[c]{ \\ with  degree 3}}  & $g_0(\mathbf{x}^*)$ & 14519.363 & 14519.363&14519.363 & 13227.616&14392.289 \\
        & $\tilde{g}_0( \mathbf{\tilde{x}})$ & 14533.311 & 14533.310&  14533.311& 13254.461& 14391.669\\
         & time(s) & 1.044&1.842 &73.693 & 603.023& 607.092\\
         & gap(\%) & 0& 0.010& 0& 41.736& -\\\midrule
        \multirow{4}{*}{\makecell[c]{ MISSOC \\ with  degree 4}} &$g_0(\mathbf{x}^*)$& 14522.985 &14522.984 & 14392.289& \multirow{5}{*}{infeas.}  & 14537.257\\
        & $\tilde{g}_0( \mathbf{\tilde{x}})$ & 14626.205& 14626.204& 14485.485 & &14623.143 \\
         & time(s) &  2.038& 3.208& 602.660& & 606.801\\
         & gap(\%) &  0.088& 0.010& 35.506& & -\\ \midrule
        \multirow{4}{*}{\makecell[c]{ MISSOC \\ with  degree 5}}  & $g_0(\mathbf{x}^*)$ &14537.257 &14537.257 & 14537.257& 5695.385& 14433.098\\
        & $\tilde{g}_0( \mathbf{\tilde{x}})$ & 14568.069 & 14568.079 & 14568.079 & 5404.909& 14433.864\\
         & time(s) & 0.823
&7.974 & 602.744 & 603.214& 607.706\\
         & gap(\%) & 0&0.010 &4.688 & 305.489& -\\ \midrule
        \multirow{4}{*}{\makecell[c]{ MISSOC \\ with  degree 6}}  & $g_0(\mathbf{x}^*)$ & 14537.257& 14537.257& 14537.098& n/a& 14392.289\\
        & $\tilde{g}_0( \mathbf{\tilde{x}})$ & 14570.065& 14569.657& 14570.167 & n/a& 14407.816\\
         & time(s) & 0.794& 5.775& 436.312& 602.953& 607.218 \\
         & gap(\%) & 0.004& 0 & 0&n/a & -\\ \midrule 
        \multirow{4}{*}{\makecell[c]{MISSOC \\ with  degree 7}}  & $g_0(\mathbf{x}^*)$ & 14537.257& 14537.25& 9835.987& \multirow{5}{*}{infeas.} & 14433.098\\
        & $\tilde{g}_0( \mathbf{\tilde{x}})$ & 14595.642&14595.650 & 8949.783 & & 14457.267\\
         & time(s) &6.669 & 117.202& 602.560& & 608.120\\
         & gap(\%) & 0.002& 0.010& 69.336& & -\\ 
         \bottomrule
\end{tabular}
\caption{Results obtained by directly solving the original Hydro Unit Commitment problem and by solving it with MISSOC using approximations of degrees 1 to 7.}
\label{tab:results_HUC}
\end{table}
}

As shown in Table~\ref{tab:results_HUC}, solving the original problem with any of the tested general-purpose MINLP solvers leads to poor results after the 600 second time limit. Gurobi declares the original problem as infeasible,  while BARON, COUENNE and BONMIN return objective values that are far below the best known feasible value of 14533.1. Recall that this is a maximization problem. These outcomes confirm the difficulty of solving the original problem and motivate addressing it with the proposed MISSOC algorithm.

Addressing  the problem through  MISSOC using the different approximations given by distinct degrees yields significantly better results in most cases. In general, the objective values obtained are much closer to that of the best known feasible solution. More importantly, the objective values in some cases even improve it. This is a very relevant result given that this instance is a real-world case study. Improving the best feasible solution known occurs with SC-MINLP and Gurobi for degrees 5 to 7, BARON for degrees 5 and 6 and BONMIN with degree 4. 

When using COUENNE in MISSOC, it does not improve the best known feasible solution for any degree. Indeed, in general, it underperforms across all degrees when compared to the other solvers. It returns the lowest (and therefore worst) objective values for degrees 2 and 3, reports infeasibility for degrees 4 and 7, fails to find a feasible solution within the time limit for degree 6 (indicated with ``n/a'') and returns a solution significantly below the best known feasible value for degree 5.

SC-MINLP and Gurobi stand out among the tested solvers in MISSOC in terms of objective values, computational times and the optimality gaps of the surrogate problems. Using these solvers in MISSOC  leads to the best objective values, except for degree 4, where using BONMIN reports a higher value. In addition, SC-MINLP and Gurobi report zero or near-zero optimality gaps, indicating that they provide stronger optimality certificates for the surrogate problem~\eqref{sec:surr_problem:tildeP}. Overall, SC-MINLP and Gurobi also yield the shortest computational times. Using them, MISSOC returns the solution in a short amount of time, with the exception of Gurobi for degree 7,  which requires more than 100 seconds. In contrast, BONMIN always reaches the time limit.  MISSOC with BARON requires moderate to long runtimes for degrees higher than 1, hitting the time limit in some cases.

Although  MISSOC with SC-MINLP and Gurobi reports similar results in terms of objective values, with zero or near-zero gaps, they differ in terms of computational time. SC-MINLP requires less time to solve the surrogate problems with degrees higher than 1. Moreover, its computational time is less sensitive to the degree used for the approximation. On the contrary, the runtime of Gurobi increases significantly with the degree, reaching over 100 seconds for the surrogate with degree 7, while SC-MINLP solves the same problem in under 7 seconds.

Comparing the results across the different   approximation degrees in MISSOC, in both SC-MINLP and Gurobi, a clear improvement is observed when the degree of the $B-$spline basis increases from 3 to 4. The objective values obtained with degree 4 are much closer to the reference value of 14533.1. For degrees above 4, this benchmark is even improved. This suggests that low-degree approximations may be insufficient for capturing the complex structure  of challenging real-world instances such as the Hydro Unit Commitment problem. Objective function approximations with higher degree might be needed to capture the nonlinearities accurately and to obtain a surrogate that better mirrors the original MINLP. This is further supported by the behavior of the values of $\tilde{g}_0( \mathbf{\tilde{x}})$,  which consistently move closer to the true objective as the degree increases.

In summary, the results highlight the practical relevance of MISSOC for addressing the Hydro Unit Commitment problem. In this case, directly solving the original problem leads to poor results. Solvers either fail to find a feasible solution or return significantly worse objective values. In contrast, MISSOC obtains high quality feasible solutions in short computational times, especially when solved with SC-MINLP, where it takes less than 7 seconds. Moreover and very importantly, MISSOC improves upon the best known feasible solution reported for this real-world instance. Altogether, these results provide strong evidence of the practical value of MISSOC for challenging real-world MINLPs.


\section{Conclusions and Perspectives}
\label{sec:conclusions}

In this work, we have introduced MISSOC, a novel algorithm to solve challenging MINLPs by building surrogate problems that are both accurate and more computationally tractable. The original complicating functions are approximated with smooth additive regression models fitted on sampled data with a $B-$spline approach. The surrogate problem is then built by replacing the original functions with their approximation through a formulation inspired by the classical multiple choice model. 

The novelty of MISSOC lies in its capability to incorporate expert or theoretical knowledge into the surrogate problem. When such information is available, our methodology allows the inclusion of constraints on the approximated function in terms of global shape or pointwise behavior. Shape constraints include bounds, monotonicity and curvature requirements imposed on the estimated function over the entire observed domain. For this, we have extended the existing methodology for shape-constrained smooth regression in the univariate case to the shape-constrained smooth additive regression model in the multivariate setting. The capability of incorporating these constraints into the estimated function makes MISSOC a data-driven and knowledge-driven approach for building surrogates of challenging MINLPs, thus filling a gap in the literature and offering a novel optimization framework for MINLPs.

The numerical results, together with the demonstrative example, show that MISSOC is a versatile approach for addressing challenging MINLPs, including problems with integer variables and where the complicating function appears either in the objective function or in the constraints. When constraints are approximated, the post processing step in MISSOC plays a key role for tackling feasibility recovery of the surrogate solution with respect to the original problem. This is illustrated in the Water Distribution Network case study.  This example shows that MISSOC can successfully handle practical MINLPs with complicating constraints in a realistic setting. The practical relevance of MISSOC is further highlighted in the Hydro Unit Commitment problem, a real-world MINLP with a complicating objective function. Beyond the computational performance, both real-world case studies also illustrate that MISSOC can build surrogates that are both data-driven and knowledge-driven.

The effectiveness of the MISSOC framework opens up several other promising research directions. Firstly, fitting a smooth additive regression model with constraints as detailed in Section~\ref{sec:shape_constrained_regression} involves estimating the weights  $\omega_j^L$ and $\omega_j^U$, $j=1,\dots,p$, which is done in this work with a data-driven approach. However, future enhancements of our method could include these weights as additional decision variables in the optimization problem~\eqref{eq:shape_constrained_additive:shape_constrained_problem}, allowing them to be tuned jointly with the regression coefficients. Furthermore, enforcing sparsity in the approximated piecewise polynomial function may further enhance the surrogate problem’s tractability, particularly in high-dimensional settings.

Another interesting direction is to provide a posteriori approximation error bounds for the approximating function. In particular, if the true function is known to be Lipschitz continuous with a certain constant, the approximating function could be estimated under the same Lipschitz constant by bounding its derivative. Thus, the approximation error would also be Lipschitz continuous and could be bounded over the whole covariate domain using the observed sample errors and the sampling density. When no prior Lipschitz information is available, the constant could instead be estimated from data, following ideas such as those in \cite{almasi5229190lipschitz}.

Moreover, extending MISSOC to black-box optimization and constraint learning would enhance its applicability in settings where analytical expressions of the objective function or constraints are not available \citep{Boukouvala2017,fajemisin2024optimization, Maragno2023mio-with-constr-learning,ten2023resource}. Finally, another worthwhile direction to explore is a dynamic approximation refinement: starting with a rough surrogate built from a small number of samples and gradually improving it in regions where good solutions are likely to be found.

\section*{Data availability statement}
The authors confirm that the data and code used to obtain the numerical results presented in Section 5 are available in public repositories linked within the article.

\section*{Acknowledgments}

This publication was supported by the Chair ``Integrated Urban Mobility'', backed by L’X - \'Ecole Polytechnique and La Fondation de l’\'Ecole Polytechnique. The Partners of the Chair shall not under any circumstances accept any liability for the content of this publication, for which the author shall be solely liable. This research benefited from the support of the FMJH Program PGMO and the grants PID2022-13724OB-I00 and CNS2023-144260, funded by MICIU/AEI /10.13039/501100011033, first also by FEDER/UE and last also by European Union NextGenerationEU/PRTR.


\end{document}